\documentclass{article}
\usepackage{amssymb}
\usepackage[babel=true,kerning=true]{microtype}

\usepackage{amsthm,leqno}

\usepackage{amsfonts,amsmath,amssymb,enumerate}
\usepackage{graphicx}
\usepackage{tikz}
\usetikzlibrary{decorations.markings}
\usetikzlibrary{plotmarks}
\usetikzlibrary{patterns}

\newcommand\res{\mathop{\hbox{\vrule height 7pt width .5pt depth 0pt
\vrule height .5pt width 6pt depth 0pt}}\nolimits}

\newcommand{\ds}{\displaystyle}

\newcommand{\be}{\begin{equation*}}
\newcommand{\ee}{\end{equation*}}
\newcommand{\beq}{\begin{equation}}
\newcommand{\eeq}{\end{equation}}
\newcommand{\begincal}{\begin{eqnarray*}}
\newcommand{\fincal}{\end{eqnarray*}}

\newtheorem{thm}{Theorem}[section]
\newtheorem{lemma}{Lemma}[section]

\newtheorem{prop}{Proposition}[section]
\newtheorem{defi}{Definition}[section]

\newtheorem{rem}{Remark}[section]

\newcommand{\Ae}{\D\setminus B(0,r)}
\newcommand{\eps}{{\varepsilon}}
\newcommand{\R}{{\mathbb R}}
\newcommand{\D}{{\mathbb D}}
\newcommand{\C}{{\mathbb C}}
\newcommand{\N}{{\mathbb N}}
\newcommand{\h}{{\mathbb H}}
\newcommand{\Z}{{\mathbb Z}}
\newcommand{\Hy}{{\mathbb H}}

\def\al{\alpha}
\def\la{\lambda}
\def\eps{\varepsilon}
\def\ds{\displaystyle}
\def\ov{\overline}

\def\p{\partial}
\def\res{\mathop{\hbox{\vrule height 7pt width .5pt 
depth 0pt\vrule height .5pt width 6pt depth 0pt}}\nolimits}

\title{Energy Quantization of Willmore surfaces at the boundary of the Moduli Space}
\author{Paul Laurain\footnote{Institut de Math\'ematiques de Jussieu,
75205 PARIS Cedex 13,
France\ } and Tristan Rivi\`ere\footnote{Department of Mathematics, ETH Zentrum,
CH-8093 Z\"urich, Switzerland.}}

\begin{document}
\maketitle
\begin{abstract}
We establish an energy quantization result for sequences of Willmore surfaces when the underlying sequence of Riemann surfaces is degenerating
in the moduli space. We notably exhibit a new residue which quantifies the potential loss of energy in collar regions. Thanks to these residues, we also establish the compactness (modulo the action of the M\"obius group of conformal transformations of ${\R}^3\cup\{\infty\}$) of the space of  Willmore immersions of any arbitrary closed 2-dimensional oriented manifold into ${\R}^3$ with uniformly bounded conformal class and energy below $12\pi$. 
\end{abstract}
\section*{Introduction}
In the study of concentration compactness phenomena  it is a central question to understand ``where and in which quantity'' the energy dissipates. The first step in analyzing such phenomena consists in looking for $\epsilon-$regularity properties. Such a property roughly says that, under a given threshold of energy, the convergence is strong. Once such a property is established, for conformally invariant problems in particular, a covering argument identifies the
points where the energy concentrates. The question is then to understand if the whole energy concentrating at these points is given exclusively  by the sum of the energies of the so called ``bubbles''  forming at these points or if there is some additional energy needed to ``anchor'' these bubbles to the rest of the solutions
in the so called ``neck region''. ``Neck regions`` are  annuli  of degenerating conformal classes separating  the bubbles between themselves or separating the bubbles with the macroscopic solution.

The Willmore energy has been introduced  in the XIX century in non linear mechanics as being the {\it ad hoc} modelization of the free
energy of a bent two dimensional elastic membrane. It was then independently introduced in geometry by Wilhelm Blaschke around 1920 in an effort to merge
minimal surfaces theory and conformal invariance. If ${\vec{\Phi}}$ denotes the immersion of an abstract closed surface  $\Sigma$ into an euclidian space
${\R}^m$, the Willmore energy of such an immersion is given by 
\[
W({\vec{\Phi}}):=\int_\Sigma|\vec{H}_{\vec{\Phi}}|^2\ dvol_{g_{\vec{\Phi}}}
\]
where $g_{\vec{\Phi}}$ is the first fundamental form of the immersion (i.e. the induced metric by ${\vec{\Phi}}$) $dvol_{g_{\vec{\Phi}}}$ is the associated volume form and $\vec{H}_{\vec{\Phi}}:=2^{-1}\,\mbox{tr}_{g_{\vec{\Phi}}}\vec{\mathbb I}_{\vec{\Phi}}$ is the half of the trace of the second fundamental form $\vec{\mathbb I}_{\vec{\Phi}}$ of ${\vec{\Phi}}$. Blaschke proved that the lagrangian $W$ is invariant under conformal transformations for closed surfaces. That is to say, for any
generic\footnote{``Generic'' means that  $\Psi^{-1}(\infty)\cap \vec{\Phi}(\Sigma)=\emptyset$.} element $\Psi$ in ${\mathcal M}({\R}^m\cup\{\infty\})$, the M\"obius group of conformal transformations of ${\R}^m\cup\{\infty\}\simeq S^m$
\[
W(\Psi\circ\vec{\Phi})=W(\vec{\Phi})\quad.
\]
 During a long time the minimal surfaces and their conformal transformations were the only known critical points of $W$. One of the reasons for such a lack of examples and progresses during almost 45 years is possibly due to the fact that the Euler-Lagrange equation of $W$ is a non-linear elliptic system of order $4$ which made it difficult to be studied from an analyst perspective at a time where this high order PDE theory was not much developed\footnote{Indeed, in a conformal parametrization $\vec{\Phi}$, the Willmore functional may be recast as
\[
W(\vec{\Phi})=\frac{1}{4}\int_{\Sigma}|\Delta_{g_{\vec{\Phi}}}\vec{\Phi}|^2\ dvol_{g_{\vec{\Phi}}}\quad,
\]
thereby giving rise to a fourth-order problem}. The seminal paper of Tom Willmore (\cite{Wi}) has relaunched the interest for the lagrangian to which his name has since then been given. 

In the present work we are interested in sequences of immersions $\vec{\Phi}_k$ of a given closed surface $\Sigma$ into ${\R}^m$ which are critical points of $W$  and
which are below a given energy level. It has been proven in \cite{Riviere08} and \cite{Riviere10} that critical points to $W$ are satisfying an $\epsilon-$regularity property (see theorem~\ref{ereg})
and that, under the assumption that the conformal class of the metrics $g_{\vec{\Phi}_k}$ is controlled, the sequence, modulo the action of the M\"obius group ${\mathcal M}({\R}^m\cup\{\infty\})$ , in a sequence of conformal parametrizations\footnote{This sequence is arbitrary as long as $\Sigma\ne S^2$ otherwise there is an {\it ad hoc} choice of gauge
in ${\mathcal M}(S^2)$ which has to be made (see \cite{MR}).},
is compact in any $C^l$ norm away from finitely many points. The iteration of this result at the various concentration points generates a bubble tree of Willmore surfaces.

In \cite{LR0} the authors established an estimate on the Green function of the Laplace operator of any degenerating sequence of constant Gauss curvature metric which permits to extend the previously mentioned concentration compactness result for Willmore surfaces to the case where the underlying conformal classes degenerate. The difficult question to know whether or not some energy could ``dissipate'' in the ``neck regions'' of the limiting ``Willmore bubble tree''
was left open in this degenerating case. The following theorem, which is the main result of the present work, asserts that, if there is some loss of energy in a collar region, the amount of this loss has an explicit expression in terms of a residue and the hyperbolic length of the collar. It can be seen as a counterpart of a result of Zhu \cite{Zhu}  for harmonic maps. Let us define this residue. Without loss of generality we can assume that  the collar region\footnote{see section 1.4 for details }  is conformally parametrized by $\D\setminus B(0,e^{-1/l_k})$ where $l_k$ is the length of shrinking geodesic which corresponds to the circle  of radius $e^{-1/2l_k}$. For $\vec{\Phi}_k :\D\setminus B(0,e^{-1/l_k}) \rightarrow \R^m$  a Willmore immersion we set , for $e^{-1/l_k}<r<1$,
\beq
2 \pi \vec{c}^{\,k}= \int_{\partial B(0,r)} \partial_\nu \vec{H}_k-3\pi_{\vec{n}_k}(\partial_\nu \vec{H}_k) - \star(\partial_\tau \vec{n}_k\wedge \vec{H}_k) \, d\sigma  
\eeq
where $\nu$  and $\tau$  are respectively  a unit normal  and a unit tangent to $\partial B(0, r)$ such that $(\nu,\tau)$ is oriented, and $\vec{n}_k$ and $\vec{H}_k$ are respectively the normal $(m-2)$-vector and the mean curvature vector of $\vec{\Phi}_k$. The operations $\star$, $\res$ and $\wedge$ are classical operations on multi-vectors that we recall in the first part of the paper. We also set

\beq
2\pi c^{\,k}_0=\int_{\gamma_k^i} - \langle (\partial_{\nu_k} \vec{H}_k-3\pi_{\vec{n}}(\partial_{\nu_k} \vec{H}_k) - \star(\partial_{\tau_k} \vec{n}\wedge \vec{H})) ,\vec{\Phi}_k\rangle  \, d\sigma
\eeq
and
\beq
2\pi \vec{c}^{\,k}_1=\int_{\gamma_k^i} - (\partial_{\nu_k} \vec{H}_k-3\pi_{\vec{n}}(\partial_{\nu_k} \vec{H}_k) - \star(\partial_{\tau_k} \vec{n}\wedge \vec{H}_k)) \wedge {\vec{\Phi}_k}  -(-1)^{m-1} 2(\star(\vec{n}_k\res \vec{H}_k))\res  \partial_{\tau_k} {\vec{\Phi}_k}\, d\sigma.
\eeq

The quantities $\vec{c}$, $c_0$ and $\vec{c}_1$ are independent of $r$ as shown in section 3.1. The quantity $\vec{c}$ already appear in \cite{Riviere08} and \cite{BR2}  as a residue which permits to erase true branch point of Willmore punctured disc( see also  \cite{KS}). But it is not conformally invariant, even not scaling invariant. For instance $\vec{c}$ is zero when compute on the closed geodesic of a catenoid, but not on the corresponding curve on the inverted catenoid, see Remark 1.1 \cite{BR2}. This remark have a very important consequence in the bubbling phenomena, see remark \ref{rm-I.1} below and theorem \ref{th-main2}. Moreover, it is  hopeless to try to find a residue, like $\vec{c}$, which measures a defect of regularity and will be conformally invariant. Taking again the example of the inverted catenoid, we can blow-up it to a union of to plane, if we get a scaling invariant residue, it must vanishes on the inverted catenoid, since the blowdo-up is smooth. Hence this residue not detect the defect of regularity\\

But $c_0$ and $\vec{c}_{1}$ are clearly invariant under the composition by isometries and dilations. They can be considered as  the Willmore analogue of  flux for CMC-surfaces (see \cite{Lo}). Since the Willmore equation is fourth order, it is not surprising to get two "fluxes". Moreover we can check, see section \ref{AHT} for details, that there exit Willmore surfaces for which those residues are non zero. Such examples are provided by considering Willmore Hopf tori, see \cite{Pinkall}. The cancellation of the residue $\vec{c}_1$ forces the generating curve (an elastica on $S^2$) to be a geodesic of $S^2$, hence the surface will be equivalent to a Clifford torus. Therefore the family of Hopf tori which are not necessary Clifford and produced by Pinkall in the above mentioned work, provides good examples. More generally, it is easy to prove that our residues vanishe on any minimal surfaces of $S^m$. Hence it provides a new way to detect non-minimal Willmore surfaces.\\

 Those quantities being defined we can state our main result.

\begin{thm}
\label{th-main}
Let $(\Sigma, h_k)$ a sequence of closed surfaces with fixed genus, constant curvature and normalized volume if needed. We assume that this sequence converges\footnote{The convergence holds in the classical sense of Mumford compactification recalled in \cite{Hummel}.} to a nodal  surface $(\widetilde{\Sigma},\widetilde{h})$ and we denote by $\{\gamma_k^i\}$ the finite number of pinching geodesics. Then let  ${\vec{\Phi}}_k : (\Sigma, h_k) \rightarrow \R^m$ a sequence of conformal Willmore immersions with bounded energy, i.e.
$$ \limsup_{k\rightarrow +\infty} W({\vec{\Phi}}_k) <+\infty$$
and such that, around every degenerating geodesic,
$$\lim_{k\rightarrow +\infty}  \frac{\vec{c}_1^{\,k}}{\sqrt{l_k}}=0.$$
 Then, denoting  $(\tilde{\Sigma}^l)_{1\leq l\leq q}$ the  connected components of $\tilde{\Sigma}$,  there exists $q$ branched smooth immersions ${\vec{\Phi}}_\infty^l :\tilde{\Sigma}^l \rightarrow \R^m$ and a finite number of possibly branched  immersions $\omega_j : S^2\rightarrow \R^m$ and $\zeta_t : S^2\rightarrow \R^m$ which are all Willmore  away from possibly finitely many points, and such that, up to a subsequence, 
\begin{equation}
\label{identity}
 \lim_{k\rightarrow +\infty} W({\vec{\Phi}}_k)=\sum_{l=1}^{q}  W({\vec{\Phi}}_\infty^l) + \sum_{j=1}^{p} W(\omega_j) + \sum_{t=1}^{q} (W(\zeta_t)-m_t 4\pi)\quad.
 \end{equation}
where $m_t$ is the integer multiplicity of $\zeta_t$ at the origin.
\hfill $\Box$
\end{thm}
\begin{rem}
\label{rm-I.1}
In \cite{BR} the second author and Bernard established the corresponding  result but under the additional assumption that the conformal class induced by the sequence
was pre-compact in the moduli space ${\mathcal M}(\Sigma)$. Under this much stronger assumptions the branched immersions $\omega_j$ and $\zeta_t$
are ``true'' Willmore surfaces in the sense that the Willmore equation is satisfied everywhere away from the branched points and moreover the first residue, $\vec{c}$, is zero on any curve surrounding these branched points. This excludes  surfaces like  the catenoid (or its inversion) in the bubble tree, but not the  Enneper surface for instance. This observation is a starting point for improving the classical bound of $8\pi$ for having compactness, see  theorem \ref{th-main2} . In theorem~\ref{th-main} instead we cannot exclude {\it a-priori} the first residue to be non zero around  the cuspidal point  for both the $\omega_j$ and the $\zeta_t$ 
\hfill $\Box$
\end{rem}

The previous remark which excludes catenoid in the  bubble trees, permits to increase the level under compactness holds true (Modulo the action of the M\"obius group). 

\begin{thm}
\label{th-main2}
Let  $\Sigma$ be a closed surface of genus $g\geq 1$  and ${\vec{\Phi}}_k : \Sigma \rightarrow \R^3$ a sequence of conformal Willmore immersions such that
$[\Phi_k^*(\xi)]$,  {the conformal class of the pullback metric, remains in a compact set of the moduli space}
and 
$$ \limsup_{k\rightarrow +\infty} W({\vec{\Phi}}_k) < 12\pi .$$
Then, There exists a diffeomorphism $\Psi_k$ of $\Sigma$ and an conformal transformation $\Theta_k$ of $\R^3\cup\{\infty\}$, such that $\Theta_k \circ \Phi_k\circ \psi_k$, up to a subsequence, converges to a smooth Willmore immersion $\Phi_{\infty}: \Sigma \rightarrow \R^3$  in $ C^\infty(\Sigma)$.
\hfill $\Box$
\end{thm}

This result was already known, when $\Sigma=S^2$, in fact there is a complete classification of Willmore sphere in $\R^3$, see \cite{Bryant}, and $\R^4$, see \cite{Montiel}. It is also known that we have compactness when $\Phi_k$ is an embedding, see \cite{Li}. But nothing, was know when $\Phi_k$ is  an immersion with energy above $8\pi$.


Coming back to theorem \ref{th-main}, the new difficulty posed by the non compactness of the underlying conformal classes in comparison with the previous quantization result in \cite{BR} will come from the formation of {\it collars}. By definition, these {\it collars}
are conformally equivalent to degenerating annuli. Unlike the {\it neck regions} from \cite{BR}, which are also conformally equivalent to degenerating annuli, 
the solution is \underbar{non extendable} throughout the interior disc of the annuli. This lack of extendability is responsible for the presence
of residues which were automatically zero in the Bernard-Rivi\`ere case.
The main achievement of the present work is to derive a control of these {\it residues} in order to prove our main results.

\medskip

The search for ``energy quantization results'' of the form above for conformally invariant problems is at the origin of numerous works in geometric analysis.
For harmonic maps ant it's heat flow for instance we can quote \cite{SaU}, \cite{St}, \cite{DiT}, \cite{Pa} or for Yang-Mills \cite{Ri5}.
 However, these problems are all of second-order elliptic or parabolic types. The novelty of \cite{BR} was to establish for the first time an energy quantization result and a full bubble-neck decomposition for a fourth-order problem. The proof of the ``energy quantization'' in \cite{BR} was using some of the   {\it integrability compensation} lemma in interpolation spaces (mostly Lorentz spaces - see the subsection 1.3 below) coming from \cite{LR} where they have been originally conceived for proving the ``energy quantization'' property for general conformally invariant second order problems in two dimensions. The use of {\it integrability by compensation} in interpolation spaces for proving ``energy quantization'' properties goes back to a work of the second author in collaboration with Fanghua Lin (see \cite{LiR}).

\medskip

In \cite{Zhu}, Miaomiao Zhu proved that for general sequences of harmonic maps from degenerating Riemann surfaces into a given manifold possessing
{\it harmonic spheres}\footnote{Harmonic spheres in a manifold $N^n$ are non constant harmonic maps from $S^2$ into $N^n$. This space is non empty for instance if $\pi_2(N^n)\ne 0$ but this condition is not necessary as the example $N^n=S^3$ shows.}  ``energy quantization'' usually \underbar{does not} hold. Therefore our result above comes as a little surprise, since we where expected two residues for a forth order problem,  but  we are able to prove that if there is some loss it can come only from the second residue $\vec{c}_1$.

\medskip

{\bf Acknowledgements} : 
{\it Part of this work has been carried out while the first author was visiting  the {\it Forschungsinstituts f\"ur Mathematik} at E.T.H. Z\"urich, he would like to thank the institute for its hospitality and the excellent working conditions.}

\section{Preliminaries}
\subsection{Notations}
Here we introduce some notations for vector calculus and differential calculus. Since we work on Riemann surfaces, it will be useful to have notation for the rotation of the gradient, so let denotes $\nabla^\bot =J\circ\nabla$ , where $J$ is the complex structure. Which gives   for $\D$,
$$\nabla^\bot f = -\frac{\partial f }{\partial y} \frac{\partial  }{\partial x} + \frac{\partial f }{\partial x} \frac{\partial  }{\partial y}$$ 

Let ${\vec{\Phi}} :\Sigma \rightarrow \R^m$ a smooth immersion and $\vec{e}$ the map which to a point $p\in \Sigma$ assigns the oriented $2$-plane given by the push forward by ${\vec{\Phi}}$ of the oriented tangent space $T_p \Sigma$. Using a positive orthonormal basis $(\vec{e}_1,\vec{e}_2)$ of ${\vec{\Phi}}_* T_p\Sigma$, we get 
$$\vec{e}=\vec{e}_1\wedge \vec{e}_2.$$
The Gauss map $\vec{n}$ assigns the oriented  $m-2$-orthogonal plane to $\vec{e}$, that is to say
$$\vec{n}=\star \vec{e} =\vec{n}_1 \wedge \dots \wedge \vec{n}_{m-2},$$
where $\star$ is the Hodge operator on $\R^m$: if $\alpha\in \bigwedge^p \R^m$ then  $\star\alpha\in \bigwedge^{m-p} \R^m$ such that for all  $\beta\in \bigwedge^{m-p} \R^m$ we get 
$$\beta \wedge \star \alpha =\langle \beta,\alpha \rangle \star 1,$$
where $\star 1$ is the canonical volume form of $\R^m$.\\

We will also need some other operator on $\bigwedge \R^m$. First, the contraction operator $\res$ : for every choice of $p-$, $q-$ and $p-q$ vectors, respectively $\alpha$, $\beta$ and $\gamma$ the following holds
$$\langle \alpha\res\beta ,\gamma\rangle = \langle \alpha,\beta\wedge\gamma\rangle.$$  
Thanks to this operator we can define the projection on the normal bundle as follow, for every $\vec{w}\in \R^m$ we set 
\beq
\label{pi}
\pi_{\vec{n}}(\vec{w})=(-1)^{m-1} \vec{n}\res(\vec{n}\res\vec{w}).
\eeq
Then this operator can be generalised assigning to a pair of $p-$ and $q-$ vectors of $\R^m$ a $p+q-2-$ vector of $\R^m$ such that, for all $\alpha\in\bigwedge^p \R^m$ and all $\beta \in\bigwedge^1 \R^m$,
\beq
\label{b1}
\alpha \bullet \beta =\alpha \res \beta 
\eeq
and, $\alpha\in\bigwedge^p \R^m$, $\beta \in\bigwedge^q \R^m$ and $\gamma \in\bigwedge^r \R^m$,
\beq
\label{b2}
\alpha \bullet (\beta \wedge \gamma)=(\alpha \bullet \beta)\wedge \gamma +(-1)^{qr}(\alpha \bullet \gamma)\wedge \beta .
\eeq

\subsection{Weak immersions}
Let $\Sigma$ a smooth compact surface equipped with a reference smooth metric $g_0$. One defines the Sobolev spaces $W^{k,p}(\Sigma,\R^m)$ of measurable maps from $\Sigma$ into $\R^m$ into the following way
$$W^{k,p}(\Sigma,\R^m)=\left\{ \vec{f}:\Sigma\rightarrow \R^m  \mathrm{ measurable\ s.t. } \sum_{l=0}^{\,k} \int_\Sigma \vert \nabla^l \vec{f}\vert^p_{g_0} \, dv_{g_0} < +\infty\right\} .$$
Since $\Sigma$ is compact it is not difficult to see that this space is independent of the choice of $g_0$ we have made.\\

Let ${\vec{\Phi}}\in W^{1,\infty}(\Sigma, \R^m)$, we define $g_{\vec{\Phi}}$ to be the following symmetric bilinear form
$$g_{\vec{\Phi}}(X,Y)= \langle d\vec{\Phi}(X), d\vec{\Phi}(Y)\rangle,$$
and we shall assume that there exists $C_{\vec{\Phi}}>1$ such that 
\beq
\label{C}
C_{\vec{\Phi}}^{-1}\ g_0(X,X) \leq g_{\vec{\Phi}}(X,X) \leq C_{\vec{\Phi}}\ g_0(X,X).
\eeq
For such a map, we can define the Gauss map as being the following measurable map in $L^\infty(\Sigma)$ taking values int the Grassmannian of oriented $m-2$-planes of $\R^m$,
$$\vec{n}_{\vec{\Phi}} = \star \frac{\frac{\partial {\vec{\Phi}}}{\partial x}\wedge \frac{\partial {\vec{\Phi}}}{\partial y}}{\left\vert \frac{\partial {\vec{\Phi}}}{\partial x}\wedge \frac{\partial {\vec{\Phi}}}{\partial y}\right\vert}.$$

We then introduce the space $\mathcal{E}_\Sigma$ of weak immersions of $\Sigma$ with bounded second fundamental form as follow:

$$ \mathcal{E}_\Sigma = \left\{\begin{array}{c}{\vec{\Phi}} \in W^{1,\infty} (\Sigma) \hbox{ which satisfies (\ref{C}) for some }C_{\vec{\Phi}}>0  \\
\\ \
\mathrm{ and } \int_\Sigma \vert  d\vec{n}_{{\vec{\Phi}}}\vert_{g_{\vec{\Phi}}}^2 \, dvol_{{\vec{\Phi}}} < +\infty \end{array}\right\},$$
where $g_{\vec{\Phi}}={\vec{\Phi}}^* \xi$.\\

It is proved in \cite{NCRiviere} that any weak immersion defines a smooth conformal structure on $\Sigma$. Let ${\vec{\Phi}}\in \mathcal{E}_\Sigma$, we denote by $\pi_{\vec{n}_{\vec{\Phi}}}$ the orthonormal projection of vector in $\R^m$ onto the $m-2$-plane given by $\vec{n}_{\vec{\Phi}}$. With these notations the second fundamental form of the immersion at $p$ is given by
$$\forall X,Y \in T_p \Sigma \ \ \vec{\mathbb{I}}_{\vec{\Phi}}(X,Y)=\pi_{\vec{n}_{\vec{\Phi}}} d^2 {\vec{\Phi}}(X,Y),$$
and the mean curvature vector of the immersion at $p$ is given by 
$$ \vec{H}_{\vec{\Phi}}=\frac{1}{2} \mathrm{tr}_{g_{\vec{\Phi}}}(\vec{\mathbb{I}}_{\vec{\Phi}}).$$
A natural quantity while considering such immersions is the Lagrangian given by the $L^2$-norm of the second fundamental form :
$$E({\vec{\Phi}})=\int_\Sigma \vert \vec{\mathbb{I}}\vert_g^2 \, dv_g.$$
An elementary computation, using Gauss-Bonnet formula, gives
$$E({\vec{\Phi}})=\int_\Sigma \vert \vec{\mathbb{I}}_{\vec{\Phi}}\vert_{g_{\vec{\Phi}}}^2 \, dvol_{g_{\vec{\Phi}}}=\int_\Sigma \vert  d\vec{n}_{{\vec{\Phi}}}\vert_{g_{\vec{\Phi}}}^2 \, dvol_{g_{\vec{\Phi}}} = 4 W({\vec{\Phi}})-4\pi \chi(\Sigma), $$
where $\chi(\Sigma)$ is the Euler characteristic and 
$$W({\vec{\Phi}})= \int_\Sigma \vert \vec{H}_{\vec{\Phi}}\vert_{g_{\vec{\Phi}}}^2 \, dvol_{g_{\vec{\Phi}}},$$
is the so called Willmore energy.\\

\subsection{Lorentz spaces}
Here we recall some classical facts about Lorentz spaces, see \cite{Gra1} for details.\\

\begin{defi}  Let $\Omega$ be a domain of $\R^{\,k}$, $p \in (1,+\infty)$ and  $q \in [1,+\infty]$. The Lorentz space $L^{p,q}(\Omega)$ is the set of measurable functions
$f : \Omega\rightarrow \R$ such that
$$\Vert f\Vert_{p,q}:= \left(\int_{0}^{+\infty} \left( t^\frac{1}{p} f^{*}(t)\right)^q \frac{dt}{t}\right)^\frac{1}{q} <+\infty \hbox{ if } q<+\infty$$
or
$$\Vert f\Vert_{p,\infty}:= \sup \left(t^\frac{1}{p} f^{*}(t)\right ) \hbox{ if } q=+\infty$$
where   $f^*$ the decreasing rearrangement of $|f|$.
\end{defi}

$\Vert\, \, \Vert_{p,q}$ happens to be a quasi norm  equivalent to a norm for which $L^{p,q}$ is a Banach space. Each $L^{p,q}$ may be seen as a deformation of $L^p$. For instance, we have
the strict inclusions
$$L^{p,1} \subset L^{p,q'} \subset L^{p,q''}\subset L^{p,\infty},$$
if $1 < q'< q''$. Moreover,
$$L^{p,p} = L^p.$$
Furthermore, if $\vert \Omega \vert$ is finite, we have that for all $q$ and $q'$,
$$p > p' \Rightarrow  L^{p,q} \subset L^{p',q'}.$$
Finally, for $p \in (1,+\infty)$ and  $q \in [1,+\infty]$, $L^{\frac{p}{p-1},\frac{q}{q-1}}$  is the dual of $L^{p,q}$.\\

{\bf Important remarks:}
Using the fact that $f^*(t)=\inf\{s>0 \hbox{ s.t. } d_f(s)\leq t\}$ where $d_f$ is the distribution function of $|f|$, we see that the $L^{2,\infty}$ norm of $f$  is finite if and only if  $\displaystyle \sup_{t>0} t^2 \left\vert \left\{ x\in D \; \vert \; \vert f (x, \, . \, ) \vert \geq t \right\}\right\vert$ is finite. Hence we easily get the following important facts:
\beq
\label{r1}
\frac{1}{\rho} \in L^{2,\infty}
\eeq
and, there exists $C>0$ such that
\beq
\label{r2}
\frac{\vert\log(r)\vert}{C} \leq \left \Vert \frac{1}{\rho} \right\Vert_{L^{2,1}(\D \setminus B(0,r))} \leq C\vert\log(r)\vert.
\eeq
We will often needs some estimates  on the mean of some functions. Thanks to (\ref{r1}) we get the following estimate. If $f$ is radial then 
\beq
\label{e3}
\int_{r}^1  \vert f\vert\; d\rho \leq \Vert f\Vert_{L^{2,1}(\D \setminus B(0,r))}  \left\Vert \frac{1}{\rho}\right\Vert_{L^{2,\infty}(\D \setminus B(0,r))} =O \left( \Vert f\Vert_{L^{2,1}(\D \setminus B(0,r))}\right) .
\eeq

\subsection{Degenerating Riemann surfaces}
Here we recall some aspects of the Deligne-Mumford's description of the loss of compactness of the conformal class for a sequence of Riemann surfaces with fixed topology, see \cite{Hummel} for details.\\

Let $(\Sigma, c_k)$ a sequence of closed Riemann surface of fixed genus $g$. If $g=0$ then the conformal class is fixed since there is only one conformal class on the sphere. If $g=1$ then, we know that, $(\Sigma, c_k)$ is conformally equivalent to 
$\ds \R^2 / \left(\frac{1}{\sqrt{\Im(v_k)}}\Z \times \frac{v_l}{\sqrt{\Im( v_k)}}\Z\right)$ where $v_k$ lies in the fundamental domain $\{z\in \C \hbox{ s.t. } \vert\Re(z)\vert\leq 1/2 \hbox{ and } \vert z \vert \geq 1\}$ of $\Hy/\mathrm{PSL}_2(\Z)$, and we say that $c_k$ degenerates if $\vert v_k\vert \rightarrow +\infty$.\\

 If $g\geq 1$, let $h_k$ the hyperbolic metric associated with $c_k$, then $(\Sigma, c_k)$ degenerates if there exits a closed geodesic whose length goes to zero. In that case, up to a subsequence, there exists
\begin{enumerate}
\item an integer $N\in \{ 1, \dots, 3g-3\}$,
\item a sequence $\mathcal{L}_k=\{ \Gamma_k^i\, ;\, i=1\dots N\}$ of finitely many pairwise disjoint simple closed geodesics of $(\Sigma,h_k)$ with length converging to zero,
\item a closed Riemann surfaces $(\overline{\Sigma},\overline{c})$,
\item a complete hyperbolic surface $(\widetilde{\Sigma},\widetilde{h})$ with $2N$ cups $\{(q^i_1,q_2^i)\, ; \, i=1\dots N\}$ such that $\widetilde{\Sigma}$ has been obtain topologically after removing the geodesic of $\mathcal{L}_k$ to $\Sigma$ and after closing each component of the boundary of $\Sigma\setminus \mathcal{L}_k$ by adding a puncture $q_l^i$ at each of these component. Moreover $\overline{\Sigma}$ is topologically equal to $\widetilde{\Sigma}$ and the complex structure defined by $\widetilde{h}$ on $\widetilde{\Sigma}\setminus \{q_l^i\}$ extends uniquely to $\overline{c}$. We can also equipped $\overline{\Sigma}$ with a metric $\overline{h}$ with constant curvature, but not necessarily hyperbolic since the genus of $\overline{\Sigma}$ can be lower than the one of $\Sigma$.

\end{enumerate}
$(\widetilde{\Sigma},\widetilde{h})$ is called the nodal surface of the converging sequence and $(\overline{\Sigma},\overline{c})$ is its renormalization. These objects are related, in the sense that, there exists a diffeomorphism $\psi_k: \widetilde{\Sigma}\setminus \{q_l^i\}  \rightarrow \Sigma\setminus \mathcal{L}_k$ such that $\widetilde{h}_k = \psi_k^* h_k$ converge in $C^{\infty}_{loc}$ topology to $\widetilde{h}$.

\subsection{Previous results : $\eps-$regularity and global control of the conformal factor.}

The first result has to do with the fact that any weak immersion with $L^2-$bounded second fundamental form defines a unique conformal structure, see  \cite{NCRiviere}.

\begin{thm}
Let ${\vec{\Phi}}$ be a weak immersion from a surface $\Sigma$ into $\R^m$ with $L^2$-bounded second fundamental form. Then there exists a constant Gauss curvature metric $h$ on $\Sigma$ and a bilipschitz homeomorphism $\Psi$ of $\Sigma$  such that ${\vec{\Phi}}\circ \Psi$ is a conformal bilipschitz immersion from $(\Sigma,h)$ into ${\R}^m$. The induced metric $g_{\vec{\Phi}}:=(\vec{\Phi}\circ\Psi)^\ast g_{{\R}^m}$ is continuous, moreover this immersion ${\vec{\Phi}}\circ\Psi$ is in $W^{2,2}(\Sigma,{\R}^m)$ and its Gauss map is in $W^{1,2}(\Sigma,Gr_{m-2}({\R}^m)$.\hfill $\Box$
\end{thm}

Assuming $\vec{\Phi}$ is conformal from the disc $\D$ into ${\R}^m$ we will denote by $\lambda$ the conformal factor, i.e.

$$e^{\lambda} =\left | \frac{\partial {\vec{\Phi}}}{\partial x_1} \right|=\left | \frac{\partial {\vec{\Phi}}}{\partial x_2} \right | ,$$
and we will denote also by $\{\vec{e}_1,\vec{e}_2\}$ the orthogonal basis of $T_{{\vec{\Phi}}(z)}\Sigma$ given by
$$ \vec{e}_i=e^{-\lambda} \frac{\partial {\vec{\Phi}}}{\partial x_i}.$$

The existence of a conformal structure is a consequence the local estimate established in \cite{Helein}.
\begin{thm}
There exists a constant $\eps_0>0$ depending only on $m$ such that for any ${\vec{\Phi}}$ weak conformal immersion from the two dimensional disc $\D$ into ${\R}^m$ satisfying
$$\int_\D \vert \nabla \vec{n}\vert^2 dz \leq \eps_0$$
then
\[
\|\lambda-\ov{\lambda}\|_{L^\infty(D(0,1/2))}\le C\ \left(  \int_\D \vert \nabla \vec{n}\vert^2+ \Vert d \lambda \Vert_{2,\infty}\right) dz,
\]
where $\ov{\lambda} =\frac{1}{\pi} \int_\D \lambda \, dz $.
\hfill $\Box$
\end{thm}

As mentioned in the introduction the starting result in the analysis of conformally invariant problems is the so called $\eps$-regularity. In the present situation this has been proved in  \cite{Riviere08} (see theorem I.5). 

\begin{thm}[$\eps$-regularity]
\label{ereg}  There exists a constant $\eps_0>0$ depending only on $m$ and for any $A>0$ a sequence of positive numbers $C_l(A)>0$ for $l\in \N^*$ 
such that for any weak  conformal immersion  ${\vec{\Phi}}: \D\rightarrow \R^m$  satisfying $$\int_\D \vert \nabla \vec{n}\vert^2 dz \leq \eps_0$$ and $\|d\lambda\|_{L^{2,\infty}({\D})}\le A$
then
\beq
\label{er1}
\Vert \nabla^l \vec{n}_{\vec{\Phi}}\Vert_{L^\infty (D(0, 1/2))} \leq C_l(A)
\left(\int_\D \vert \nabla \vec{n}_{\vec{\Phi}}\vert^2 \, dz \right)^\frac{1}{2}.
\eeq
\hfill $\Box$
\end{thm}
The following result is a consequence of the $\eps-$regularity
\begin{thm}[Theorem I.5 of \cite{Riviere08}]
\label{conv}  There exists a constant $\eps_0>0$ and some constant $C_m$ for $m\in \N$, such that, if ${\vec{\Phi}}_k$ is a sequence of conformal Willmore  immersion of $\D$ into $\R^m$ which satisfies
$$\int_\D \vert \nabla \vec{n}_k\vert^2 dz \leq \eps_0\quad\mbox{ and }\quad\limsup_{k\rightarrow +\infty}\|d\la\|_{L^{2,\infty}({\D})}<+\infty$$
where $\vec{n}_k$ is normal associated to ${\vec{\Phi}}_k$,  then up to a dilation in the image,  ${\vec{\Phi}}_k$ converges in $C^{2}_{loc}(\D)$ to a conformal Willmore immersion ${\vec{\Phi}}_\infty$.
\hfill $\Box$\end{thm}

An other ingredient consists in controlling the conformal factor independently of the conformal class when the $L^2$-norm of the second fundamental form is bounded. This correspond to theorem 3.1 of Laurain-Rivi\`ere \cite{LR0}.

\begin{thm}
\label{a1}
 Let $(\Sigma,c_k)$ be a sequence of closed Riemann surface of fixed genus greater than one. Let denote $h_k$ the metric with constant curvature \footnote{ equal to $-1$,$0$ or $1$ and volume equal to one in the torus case} in $c_k$ and ${\vec{\Phi}}_k$ a sequence of weak conformal immersion of $\Sigma$ into $\R^m$, i.e.
$${\vec{\Phi}}_k^\ast g_{{\R}^m} =e^{2u_k}h_k,$$
where $u_k\in L^\infty(\Sigma)$. Then there exists a finite conformal atlas $(U_i,\psi_i)$ independent of $k$ and a positive constant $C$ depending only on the genus of $\Sigma$, such that 
$$\Vert d \lambda^i_k \Vert_{L^{2,\infty}(V_i)}\leq C\ W({\vec{\Phi}}_k),$$
where $\lambda^i_k$ is the conformal factor of ${\vec{\Phi}}^{\,k} \circ \psi_i^{-1}$ in $V_i=\psi_i(U_i)$, i.e. $\lambda_k^i=\frac{1}{2} \log \left\vert \frac{\partial {\vec{\Phi}}^{\,k}\circ \psi_i^{-1}}{\partial x} \right\vert=\frac{1}{2} \log \left\vert \frac{\partial {\vec{\Phi}}^{\,k}\circ \psi_i^{-1}}{\partial y}\right\vert$.
\hfill $\Box$
\end{thm}

\section{Detecting Bubbles, Necks and Collars.}
\label{preli}
We consider a sequence ${\vec{\Phi}}_k$ in ${\mathcal E}_\Sigma$  critical points of $W$ and bounded energy, i.e.
$$\limsup_k W({\vec{\Phi}}_k) <+\infty .$$

Let denote $c_k$ the conformal class defined by ${\vec{\Phi}}_k$. Let $h_k $ be a constant curvature metric such that $g_{\vec{\Phi}_k}= e^{\alpha_k} h_k$. Moreover in the case of genus one, we normalize the area of $(\Sigma, h_k)$ to be $1$.\\

In the sphere case there is only one conformal class,  we can directly apply the main result of Bernard and Rivi\`ere \cite{BR}.\\

So we don't consider the case $\Sigma=S^2$ and we shall decompose our surface into thin and thick parts. On the thick parts the metric converges smoothly and the classical theory of \cite{BR} applies. The thin parts are conformally equivalent to  long cylinders. Then we need  an equivalent result to the classical bubble tree decomposition  in this context.\\

\subsection{Bubble Tree lemma}
In this section, we generalized the bubble extraction made in section III of \cite{BR} to a collar region.
\begin{lemma}
\label{btl} Let $l_k\rightarrow +\infty$ and  ${\vec{\Phi}}_k : S^1\times [0,l_k] \rightarrow \R^m$ a sequence of Willmore immersions with $L^2$-bounded second fundamental form.
We assume that there is no concentration at the boundary that is to say for every $R>0$, ${\vec{\Phi}}_k$  and   ${\vec{\Phi}}_k( \theta, l_k-t)$ converge in $C^2(S^1\times [0,R])$. Then, either
$$\lim_{R\rightarrow +\infty }\lim_{n\rightarrow +\infty} \sup_{t\in [R,l_k-R]} \int_{S^1 \times [t,t+1]} \vert \nabla \vec{n}_k\vert^2 \, d\theta dt =0$$
or there exist $p>0$, $2p$ sequences of real $0a_k^1\leq b_k^1\, a_k^2\leq b_k^2, \dots, a_k^p \leq b_k^p$, such that 
\begin{itemize}
\item $\displaystyle \lim_{k \rightarrow +\infty } b_k^i-a_k^i > 1$ for all $1\leq i\leq p$,
\item $\displaystyle \lim_{k \rightarrow +\infty } a_k^{i+1}-b_k^i =+\infty$ or all $1\leq i\leq p-1$,

\item $\displaystyle  \lim_{k \rightarrow +\infty } \frac{b_k^i-a_k^i}{l_k} = 0 $

\item  $\displaystyle \lim_{R\rightarrow +\infty }\lim_{n\rightarrow +\infty} \sup_{t\in [b_k^{i}+R,a_k^{i+1}-R]} \int_{S^1 \times [t,t+1]} \vert \nabla \vec{n}_k\vert^2 \, d\theta dt =0 \hbox{ for all } 1\leq i\leq p-1$
\end{itemize}
and $ {\vec{\Phi}}_k^i(\theta,t) ={\vec{\Phi}}_k\left( \theta,\frac{a_k^i+b_k^i}{2}+t \right) $ satisfies a non trivial\footnote{In the sens that there is at least one bubble.} energy identity\footnote{In the sens of the main theorem of \cite{BR}.} on $\left[ \frac{a_k^i-b_k^i}{2},\frac{b_k^i-a_k^i}{2} \right]$.
\hfill $\Box$
\end{lemma}
\noindent{\bf Proof of lemma \ref{btl}:}\\

Let $\displaystyle \Gamma = \lim_{R \rightarrow +\infty } \lim_{k \rightarrow +\infty } \sup_{t\in [R,l_k-R]}\int_{S^1 \times  [t,t+1]} \vert \nabla \vec{n}_k\vert^2 \, dz$. Either $\Gamma=0$ and there is nothing to extract or $\Gamma>0$. Hence we pick $t_k\in (0,l_k)$ such that 

$$\displaystyle \lim_{k\rightarrow +\infty} t_k= \lim_{k\rightarrow +\infty} l_k-t_k=+\infty$$
and
$$ \lim_{k \rightarrow +\infty } \int_{S^1 \times  [t_k,t_k+1]} \vert \nabla \vec{n}_k\vert^2 \, dz =\Gamma >0 .$$

Then we consider  $ {\vec{\Phi}}_k^1(\theta,t) ={\vec{\Phi}}_k\left( \theta,t_k+t \right) $. Since the energy is finite, there is a finite number of concentration points, where the energy identity is satisfied since locally we can apply the main result of Bernard Rivi\`ere. So all the points of concentration of ${\vec{\Phi}}_k^1$ are contained in $[-R_1,R_1]$ for some $R_1>0$ and ${\vec{\Phi}}_k^1$ converge to some  bubble in $C^2_{loc}((S^1\times \R)\setminus \{\hbox{concentration points} \})$ .\\

Hence, setting $a_k^1=t_k-R_1$ and $b_k^1=t_k+R_1$,  the boundary hypothesis of the lemma on $S^1\times [0,a_k^1]$ and $S^1\times [b_k^1, l_k]$ are satisfied and we have an energy identity on $S^1 \times [a_k^1,b_k^1]$. Then we apply the process recursively to $S^1\times [0,a_k^1]$ and $S^1\times [b_k^1, l_k]$. The process has to stop, since at each step the central cylinder contains at least $\eps_0$ of energy since either there is concentration or one converges in $C^2_{loc}(S^1\times \R$)) to a nontrivial bubble. \hfill$\square$

\subsection{Choosing the thin part}

\subsubsection{The torus case}

The theorem 0.2 of \cite{LR0} ensures that a torus, which is isometric to the cylinder $C_k=\frac{1}{\sqrt{2\pi l_k}} \left( S^1 \times \left[0,l_k\right]\right)$ with the standard identification of its boundary components, admits the following chart 
$$
\begin{array}{cccc}
\psi_k :& \D\setminus B(0,e^{-l_k}) & \rightarrow & C_l \\
& (\theta, r) &\mapsto & \left(\frac{cos(\theta)}{\sqrt{2\pi l_k}}, \frac{sin(\theta)}{\sqrt{2\pi l_k}}, \frac{-\log(r)}{\sqrt{2\pi l_k}}\right) \, ,
\end{array}
$$
such that, the conformal factor $u_k$ of ${\vec{\Phi}}_k={\vec{\Phi}}_k \circ \psi_k$, i.e. ${\vec{\Phi}}_k^*(\xi)=e^{2u_k}dz^2$, satisfies
$$\Vert \nabla u_k \Vert_{L^{2,\infty}} \leq C.$$

Moreover we can choose the place where we "cut" the torus into a cylinder in a way that there is no concentration near the boundary.\\

Indeed there is a finite number of $t^i\in \left[-\frac{l_k}{2},\frac{l_k}{2}\right]$ such that $\ds \lim_{k\rightarrow +\infty} \int_{S^1\times [t^i,t^i+1]} \vert \nabla \vec{n}_k\vert ^2 \, dtd\theta \geq \frac{\eps_0}{2}$, then we pick $t_k$ such that $\ds \lim_{k \rightarrow +\infty} \vert t_k -t^i\vert =+\infty$. Hence, setting $\tilde{{\vec{\Phi}}}_k={\vec{\Phi}}_k(\, .\, +t_k)$, thanks to $\eps$-regularity, we have the convergence of $\tilde{{\vec{\Phi}}}_k$ in $C^2_{loc}$ to some (possibly trivial) bubble. Hence cutting the torus at $t^{\,k}$ provide a cylinder with no concentration near the boundary. 

\subsubsection{The hyperbolic case}

Thanks to the collar lemma we know that choosing $\delta < \sinh(1)$, the thin part  $\{x \in (\Sigma, h_k)\, \vert \, inj(x) <\delta\}$ consists of a finite number of collars. Up to extraction, this number is fixed for $k$ large enough, but in order to simplify notations in the rest of this part we will assume that there is only one collar. This collar contains  a smallest closed geodesic  and  is conformal to an hyperbolic cylinder of the form

$$
A_l=\left\{ z=re^{i{\vec{\Phi}}}\in \h : 1\leq r \leq e^l, arctan\left(sinh\left(\frac{l}{2}\right)\right)<{\vec{\Phi}}<\pi - arctan\left(sinh\left(\frac{l}{2}\right)\right)\right\},
$$
where the geodesic correspond to $\left\{ r e^{i\frac{\pi}{2}}\in \h : 1 \leq r \leq e^l\right\}$ and the line $\{r=1\}$ and $\{ r=e^l\}$ are identified via $z\mapsto e^l z$. This is the {\bf collar region}. It is sometimes easier to consider the following  cylindrical parametrization, i.e.

$$P_l=\left\{ (t,\theta) : \frac{2\pi}{l}arctan\left(sinh\left(\frac{l}{2}\right)\right)<t<\frac{2\pi}{l} \left(\pi - arctan\left(sinh\left(\frac{l}{2}\right)\right)\right), 0\leq \theta \leq 2\pi\right\}$$
in this parametrization the constant scalar curvature metric reads
$$
ds^2=\left(\frac{l}{2\pi sin(\frac{lt}{2\pi}) }\right)^2 (dt^2 +d\theta^2),
$$
where the geodesic corresponds to $\{ t=\frac{\pi^2}{l}\}$ and the line $\{\theta=0\}$ and $\{ \theta=2\pi\}$ are identified.

Then, as $l_k$, the length of the degenerating geodesic, goes to zero, $P_{l_k}=[0,T_k]\times S^1$ becomes a long cylinder.

Let $\psi_k$ be the chart from the cylinder to the collar and keep denoting ${\vec{\Phi}}_k ={\vec{\Phi}}_k \circ \psi_k$. Choosing $\delta$ small enough, we can assume that there is no concentration near the boundary of the collar, i.e., for  every $R>0$, ${\vec{\Phi}}_k$  and  ${\vec{\Phi}}_k(T_k-t,\theta)$ converges in $C^{2}_{loc}([0,R]\times S^1)$.

\subsection{Extraction of necks inside the collars}

From the last section, we reduce the study of the thin part to a long cylinder without concentration near the boundary. A priori, the metric is not flat, but since the energy is invariant by conformal change of metric, we can apply directly the bubble tree lemma.Indeed, lemma \ref{btl} permits us to split our cylinder in parts where we have energy identities, the $S^1\times [a^i_k, b_k^{i}]$,  and necks, the  $S^1\times [b^i_k, a_k^{i+1}]$. Then in the rest of the paper we will concentrate on a specific neck region. In order to prove that there is no energy in such a neck we will use the $L^{2,\infty}-L^{2,1}$ duality. As remark in Laurain-Rivi\`ere \cite{LR0}, this norm are no more invariant by conformal change (beside  simple dilatations and isometries of course). So we are going to fix the chart once for all. But thanks to Laurain-Rivi\`ere \cite{LR0} we know that there is an appropriate choice to make use of this theory. We make it precise in the next section.

\subsection{ $L^{2,\infty}$ estimate in the neck }
\label{2.4}
Thanks to theorem \ref{a1}, we know that we can choose our chart such that the  $L^{2,\infty}$ norm of the gradient of the conformal factor is uniformly bounded. Moreover, every collar (resp. long and thin cylinders) has a chart of annular type, i.e.
$$ A_k=\D\setminus B(0,e^{-l_k}) , $$
where $l_k \rightarrow +\infty$. In this setting the result corresponding to lemma VII.1 of \cite{BR} holds. Precisely we have the following lemma

\begin{lemma}
\label{leml2i}
Let $m\geq 3$ and  ${\vec{\Phi}}_k: \D\setminus B(0,\eps_k) \rightarrow \R^m$ be a sequence of  conformal Willmore immersion, with $\eps_k \rightarrow 0$, satisfying

$$\Vert \nabla \lambda_k \Vert_{L^{2,\infty}(\D\setminus B(0,\eps_k) )} + \Vert \nabla \vec{n}_{{\vec{\Phi}}_k} \Vert_{L^{2}(\D\setminus B(0,\eps_k) )} \leq C,$$
where $\lambda_k$ is the conformal factor, and such that 
\beq
\label{nl2i}
\lim_{R\rightarrow +\infty }  \lim_{k\rightarrow +\infty }\sup_{R \eps_k <\rho< \frac{1}{R}} \int_{B(0,2\rho) \setminus B(0,\rho)} \vert \nabla \vec{n}_{{\vec{\Phi}}_k} \vert^2 \, dz =0.
\eeq
Then for every $\eps>0$, there exists $R>0$, depending only on $C$ and $\eps$,  such that 
$$   \lim_{k\rightarrow +\infty }\sup_{R \eps_k <\vert z \vert < \frac{1}{R}} \vert z\vert   \vert \nabla  \vec{n}_{{\vec{\Phi}}_k} \vert < \eps.$$
in particular
\beq
\lim_{R\rightarrow +\infty }  \lim_{k\rightarrow +\infty } \Vert  \nabla  \vec{n}_{{\vec{\Phi}}_k} \Vert_{L^{2,\infty}(B(0,\frac{1}{R})\setminus B(0,R\eps_k) )} =0 .
\eeq
\hfill $\Box$
\end{lemma}
\subsection{ Behaviour of the conformal factor in a neck region}

First of all, let us remind  the following lemma which permits to control the behaviour of the conformal factor in a annulus.

\begin{lemma}[lemma V.3 of \cite{BR}]
\label{lambde} 
There exists $\eps_0>0$ with the following property. Let $0<r<\frac{1}{4}$. If ${\vec{\Phi}}$ is any conformal weak immersion of $\D\setminus B(0,r)$ into $\R^m$ with $L^2$-bounded second fundamental form, and satisfying
$$\Vert \nabla \vec{n}\Vert_{L^{2,\infty}} \leq \eps_0,$$
then there exist $d,A\in \R$ such that 
$$\Vert \lambda(x)-dLog(\vert x\vert)  -A\Vert_{L^{\infty}(B(0,\frac{1}{2})\setminus B(0,2r))} \leq C\left( \Vert \nabla \lambda \Vert_{L^{2,\infty}} +\int_{\D\setminus B(0,r)} \vert \nabla \vec{n}\vert^2 \, dz\right) ,$$
where $d$ satisfies
$$\left\vert 2\pi d-\int_{\partial B_r} \frac{\partial \lambda}{\partial \rho} \, d\sigma \right\vert \leq C \left( \int_{B(0,2r)\setminus B(0,r)} \vert \nabla \vec{n}\vert^2 \, dz + \frac{1}{Log(1/r)} \left( \Vert  \nabla \lambda\Vert_{2,\infty} + \int_{\D\setminus B(0,r)} \vert \nabla \vec{n}\vert^2 \, dz\right)\right),$$ 
where $C$ depends only on $m$ and $\lambda$ is as in lemma \ref{XX}.  
\hfill $\Box$
\end{lemma}

As we will see later, we need to prevent $d$ to be close to $-1$, which correspond to the parametrization of a long thin cylinder. Hence in the en of this paragraph we exclude the possibility for a sequence of $\vec{\Phi}_k$ satisfying the hypothesis of lemma \ref{lambde} and (\ref{nl2i}) to have its corresponding $d_k$ that converge to $-1$.\\

Let $\vec{\Phi}_k$ satisfying the hypothesis of lemma \ref{leml2i}.  We set $\vec{\Psi}_k= \frac{\vec{\Phi}_k(\sqrt{\eps_k} z) -c_k}{L_k}$ where  $c_k$ and $L_k$ are respectively the center of mass and  the length of $\theta \mapsto \vec{\Phi}_k (\sqrt{\eps_k} e^{i\theta})$.\\

Let $K$ a compact set of $\C\setminus\{0\}$, for $k$ large enough, $\vec{\Psi}_k$ is well defined on $K$ and its normal satisfies, thanks to  lemma \ref{leml2i}, the following estimate

\be
\label{normpsi}
 \vert \nabla \vec{n}_{\vec{\Psi}_k} \vert = o(1).
 \ee

Moreover, since  the energy does not concentrate our new conformal factor satisfies an Harnarck inequality, see Lemma V.2 of \cite{BR}. Using lemma \ref{lambde}, we can write it as $\lambda_{\vec{\Psi}_k}  = d_k \ln(\rho) +A_k+B_k$, where $d_k$ and $A_k$ are some real constant and $B_k$ is a uniformly bounded function. Since, $\vec{\Psi}_k$ is conformal and   the length of the image of $S^1$ is one, then $A_k$ is uniformly bounded. Moreover, $d_k$ is also uniformly bounded by Harnarck inequality.\\

Finally $\vec{\Psi}_k$, up to an extraction , converge to $\vec{\Psi}_\infty$ a conformal parametrization of a piece of plane in $C^2_{loc}(\C\setminus\{0\})$, and it is conformal factor satisfies
\beq
\label{di}
\frac{1}{C} \rho^{d_\infty} \leq e^\lambda \leq  C \rho^{d_\infty},
\eeq
for some $C>0$.\\

Up to orientation, we can consider $ \vec{\Psi}_\infty$ as an holomorphic  map from $\C\setminus\{0\}$ to $\C$, we note it simply $\psi$ from now. Then $\psi' dz$ is an holomorphic one form on $\C\setminus \{0\}$. Thanks to (\ref{di}),  it is in fact rational form on $\hat{\C}$. Then we get  that $d_\infty\in \Z$ and, up to multiply $\psi$  by a non null complex number, $\psi'(z)= z^{d_\infty}$. Of course the case $d_\infty=-1$ is exclude since $\frac{dz}{z}$ has no primitive on $\C\setminus\{0\}$. \\

\begin{rem}
\label{remd}  When considering a sequence of Willmore immersion, in all neck we have the decomposition of the conformal factor as in lemma \ref{leml2i}, moreover we can assume that $\vert d_k+1\vert >\frac{1}{2}$. 
\end{rem}

\section{Preliminary Estimates in Collar Regions.}
We shall very often omit to write explicitly the subscript $k$ when there is no ambiguity for understanding the argument.  

\subsection{ $L^{2,1}$-estimate on degenerating annuli}

Thanks to the conclusion of the section \ref{2.4}, in order to prove the quantification of the energy in the collar it suffices to prove an $L^{2,1}$-estimate on $\nabla \vec{n}$ in degenerating annuli.  This is however not exactly what we shall establish. Here is the main difference with case  where the conformal class stays bounded : our annuli are parts of the collars and cannot be filled inside by the solution which extends to the whole disc.  Bernard and the second author where facing a somehow similar situation in their study of Willmore surfaces near branched points (see \cite{BR2}). In order to overcome the impossibility to extend the solution throughout the middle region at each steps of the derivation of the conservation laws issued from the application of Noether theorem (see \cite{NCRiviere}) residues will show up.
 The main task of the present work will be to control these residues.\\

Let ${\vec{\Phi}} : \Ae \rightarrow \R^m$ a conformal Willmore immersion. Thanks to theorem 1.1 of \cite{Riviere08}, denoting

\beq
\label{defc}
2 \pi \vec{c}
 :=\int_{\partial \D}  \vec{\nu}.\left(\nabla \vec{H}-3\pi_{\vec{n}}(\nabla \vec{H}) + \star(\nabla^\bot \vec{n}\wedge \vec{H})\right)\, d\sigma  \eeq
there exists $\vec{L}: \Ae \rightarrow \R^m$ such that 
\beq
\label{eqL}
\nabla^\bot \vec{L} = \nabla \vec{H}-3\pi_{\vec{n}}(\nabla \vec{H}) + \star(\nabla^\bot \vec{n}\wedge \vec{H}) -\vec{c} \,\nabla  \log(\rho)
\eeq 
and following \cite{Riviere08} we observe that
\be
\left\{
\begin{split}
&\displaystyle\mbox{div}(\langle \vec{L}, \nabla^\bot {\vec{\Phi}} \rangle- \langle \vec{c} , {\vec{\Phi}}\rangle \nabla \log(\rho) )=0\quad,\\[5mm]
&\displaystyle\mbox{div}(\vec{L}\wedge \nabla^\bot {\vec{\Phi}} +(-1)^{m-1}2(\star(\vec{n}\res \vec{H})\res \nabla^\bot {\vec{\Phi}} -\vec{c} \wedge {\vec{\Phi}} \, \nabla \log(\rho)))=0\quad.
\end{split}
\right.
\ee
Let
\beq
\label{c0}
 2 \pi c_0:= -\int_{\partial \D}  \langle \vec{L}, \partial_{\tau}  {\vec{\Phi}} \rangle\, d\sigma  -\int_{\partial \D}  \langle \vec{c},   {\vec{\Phi}} \rangle\, d\sigma 
\eeq
and
\beq
\label{c1}
2 \pi \vec{c}_1:= \int_{\partial \D} \left(-  \vec{L}\wedge \partial_{\tau} {\vec{\Phi}} -(-1)^{m-1}2(\star(\vec{n}\res \vec{H})\res \partial_\tau {\vec{\Phi}}) \right) \, d\sigma -\int_{\partial \D} \vec{c}\wedge{\vec{\Phi}}  \, d\sigma  ,
\eeq
where $\partial_\tau=\frac{1}{\rho}\frac{\partial}{\partial \theta}$ and $d\sigma=\rho\, d\theta$.

Then, thanks to Poincar\'e lemma,  there exists $S :\Ae \rightarrow \R$ and $\vec{R} :\Ae \rightarrow \R^m$ such that 
\beq
\label{e1}
\nabla^\bot S = \langle \vec{L}, \nabla^\bot {\vec{\Phi}}\rangle - C_0 \nabla\log(\rho) 
\eeq
and

\beq
\label{e2}
\nabla^\bot \vec{R} =  \vec{L}\wedge \nabla^\bot {\vec{\Phi}} +(-1)^{m-1}2(\star(\vec{n}\res \vec{H})\res \nabla^\bot {\vec{\Phi}})-\vec{C}_1 \nabla \log(\rho).
\eeq
where
$C_0 =c_0 +\langle \vec{c}  ,{\vec{\Phi}}\rangle$ and $\vec{C}_1 =\vec{c}_1 +\vec{c}\wedge {\vec{\Phi}} $.

The rest of this section consists essentially in estimating the {\it residues} $C_0$ and $\vec{C}_1$ independently of the conformal class of the annuli.  This is the novelty of this paper,  those {\it residues} do not appear when the neck regions is not part of collar  since the solution is extendable throughout the internal small disc. Before doing so, we generalize the estimate (IV.24) of \cite{BR}, to our setting. We propose in fact here a shorter argument than the one in \cite{BR}.\\

\subsubsection{Pointwise estimate of $\vec{L}$}
We are  going to use the  lemma \ref{lambde} to derive a pointwise control of ${\vec{L}}$. Precisely we have.
\begin{lemma}
\label{lemma-linftyestim}There exists $C>0$, independent of $r$, such that
\beq
\label{eL}
e^{\lambda_\rho} \Vert  \vec{L} \Vert_{L^{\infty} (\partial B(0,\rho))} \leq \frac{C}{\rho} \quad  \hbox{ for all } \rho \in (2r,1/2),
\eeq
where $\displaystyle \lambda_\rho =\sup_{\vert x\vert= \rho} \lambda(x)$. \hfill $\Box$
\end{lemma}

\noindent{\bf Proof of lemma~\ref{lemma-linftyestim}.}

First of all, let us remind the pointwise estimate on $\nabla \vec{H}$. From (VI.12) of \cite{BR} and the Harnack inequality satisfies by the conformal factor, see lemma V.2 of \cite{BR}, we get that
\beq
\label{er2}
\vert \nabla \vec{H}(x) \vert \leq C e^{-\lambda_{\vert x\vert }} \frac{1}{ \vert x\vert^2}\left( \int_{B(0,2\vert x\vert)\setminus B(0,\frac{\vert x\vert}{2})} \vert \nabla \vec{n}\vert^2 \, dz \right)^\frac{1}{2} \hbox{ for all } x\in B(0,1/2)\setminus B(0,2r),
\eeq
where $C$ depends only on $\Vert \nabla \lambda\Vert_{2,\infty}$.Then we set $\vec{L}_\rho =\frac{1}{2\pi} \int_0^{2\pi} \vec{L}(\rho,\theta) \, d\theta$. Thanks to (\ref{defc}), (\ref{eqL}), (\ref{er1}) and (\ref{er2}), we get, for all $\rho\in (2r,1/2)$, 

$$\vert \vec{L}-\vec{L}_\rho\vert \leq  \int_{0}^{2\pi} \vert \nabla \vec{L}\vert \,\rho d\theta \leq C \frac{e^{-\lambda_\rho}}{\rho} \left( \int_{B(0,2\rho)\setminus B(0,\frac{\rho}{2})} \vert \nabla \vec{n}\vert^2 \, dz \right)^\frac{1}{2}.$$
Then we estimate $\vec{L}_\rho$. First of all, Thanks to (\ref{eL}), (\ref{er1}) and (\ref{er2}) we have
$$\left\vert \frac{d \vec{L}_\rho}{d\rho} \right\vert \leq C( \vert \nabla \vec{H}\vert + \vert \nabla \vec{n} \vert \vert \vec{H}\vert) \leq C \frac{e^{-\lambda_\rho}}{\rho^2}  \left( \int_{B(0,2\rho)\setminus B(0,\frac{\rho}{2})} \vert \nabla \vec{n}\vert^2 \, dz \right)^\frac{1}{2}. $$
 But, since $\vec{L}$ is defined up to a constant, we can assume that $\vec{L}_{2r}=0$ (or $\vec{L}_{1/2}=0$, this will be make clear later), which gives
\be 
\begin{split}
\vert \vec{L}_\rho \vert &\leq C  \int_{2r}^\rho \frac{e^{-\lambda_t}}{t^2}  \left( \int_{B(0,2t)\setminus B(0,\frac{t}{2})} \vert \nabla \vec{n}\vert^2 \, dz \right)^\frac{1}{2}\, dt \\
&\leq C e^{-\lambda_\rho}  \int_{2r}^\rho \frac{e^{\lambda_\rho-\lambda_t}}{t^2}  \left( \int_{B(0,2t)\setminus B(0,\frac{t}{2})} \vert \nabla \vec{n}\vert^2 \, dz \right)^\frac{1}{2}\, dt. 
 \end{split}
 \ee
 Thanks to lemma \ref{lambde} we have $\lambda_\rho-\lambda_t =d\log\left(\frac{\rho}{t}\right) +B$, where $B$ is uniformly bounded, hence, if  $d+1\leq0$,
 
\be 
\begin{split}
\vert \vec{L}_\rho \vert &\leq C e^{-\lambda_\rho}  \int_{2r}^\rho \frac{1}{t^2} \left(\frac{\rho}{t}\right)^d  \left( \int_{B(0,2t)\setminus B(0,\frac{t}{2})} \vert \nabla \vec{n}\vert^2 \, dz \right)^\frac{1}{2}\, dt. \\
 &\leq C \frac{e^{-\lambda_\rho}}{\rho}  \left(\frac{\rho}{2r}\right)^{d+1}   \int_{2r}^\rho \frac{1}{t^2}  \left( \int_{B(0,2t)\setminus B(0,\frac{t}{2})} \vert \nabla \vec{n}\vert^2 \, dz \right)^\frac{1}{2}\, tdt \\
&\leq C \frac{e^{-\lambda_\rho}}{\rho}  \left(\frac{\rho}{2r}\right)^{d+1} 
 \end{split}
 \ee
which gives the desired estimate if $d+1\leq0$. Then if $d+1\geq 0$ it suffies to assume that $\vec{L}_{\frac{1}{2}}=0$ and to integrate between $\rho$ and $1/2$.\hfill$\square$ 
 
\subsubsection{Estimating the first residue  $C_0$.}

\begin{lemma}
\label{c0est}
There exists $C>0$, independent of $r$, such that
\beq
\label{finalc0}
\left\Vert \frac{C_0}{\rho} \right\Vert_{L^{2,1}(B(0,\frac{1}{2})\setminus B(0,2r)} \leq C.
\eeq
\hfill $\Box$
\end{lemma}

\noindent{\bf Proof of lemma~\ref{c0est}.}
Integrating by part (\ref{c0}) and using Stokes's theorem, we get,  for every $z \in B(0, \frac{1}{2})\setminus B(0,2r)$, 
$$2\pi C_0(z) =\int_{\partial B(0, \vert z\vert)}   \langle\partial_{\tau} \vec{L},   {\vec{\Phi}} \rangle \, d\sigma -\int_{\partial B(0, \vert z\vert)}   \left\langle\frac{\vec{c}}{\vert z\vert},   {\vec{\Phi}} \right\rangle \, d\sigma+2\pi \langle \vec{c}  ,{\vec{\Phi}}\rangle $$
Thanks to (\ref{eqL}), we have the fact that $\partial_\tau \vec{L} =-\vec{T}_\nu +\frac{\vec{c}}{\vert z\vert}$ where $\vec{T}=\nabla \vec{H}-3\pi_{\vec{n}}(\nabla \vec{H}) + \star(\nabla^\bot \vec{n}\wedge \vec{H})$. Then
\beq
\label{c00}
\begin{split}
2\pi C_0(z) &=\int_{\partial B(0, \vert z\vert)}   -\langle \vec{T}_{\nu},   {\vec{\Phi}} \rangle \, d\sigma   + \int_{\partial B(0, \vert z\vert)}   \left\langle\frac{\vec{c}}{\vert z\vert},   {\vec{\Phi}} \right\rangle \, d\sigma  -\int_{\partial B(0, \vert z\vert)}   \left\langle\frac{\vec{c}}{\vert z\vert},   {\vec{\Phi}} \right\rangle \, d\sigma+2\pi \langle \vec{c} ,{\vec{\Phi}}\rangle \\
&=\int_{\partial B(0, \vert z\vert)}   \langle \vec{T}_{\nu}(y),   {\vec{\Phi}}(z)-{\vec{\Phi}}(y) \rangle \, d\sigma 
\end{split}
\eeq

On the one hand
\beq
\label{c01}
 \vec{T}_\nu=   \frac{\partial \vec{H}}{\partial\rho } -3\pi_{\vec{n}}\left( \frac{\partial \vec{H}}{\partial\rho }\right) - \star\left(\frac{\partial \vec{n}}{\partial \tau}\wedge \vec{H}\right)=  \pi_{T}\left(\frac{\partial \vec{H}}{\partial\rho }\right) -2\pi_{\vec{n}}\left( \frac{\partial \vec{H}}{\partial\rho }\right) - \star\left(\frac{\partial \vec{n}}{\partial \tau}\wedge \vec{H}\right),
\eeq
where $\pi_{T}$ is the projection onto the tangent plane. Writing, locally, $\ds \vec{H}=\sum_{\alpha=1}^{m-2} H_\alpha \vec{n}^\alpha$, we easily get that 
\beq
\label{pit}
\left\vert \pi_{T}\left(\frac{\partial \vec{H}}{\partial\rho }\right) \right\vert \leq \vert H\vert \vert \nabla \vec{n}\vert.
\eeq
On the other hand, let $\Gamma_{y,z}$ the direct circular arc from $y$ to $z$, we get 
\be
\begin{split}
\left\vert\left\langle \pi_{\vec{n}}\left( \frac{\partial \vec{H}}{\partial\rho }\right) (y) , {\vec{\Phi}}(z)-{\vec{\Phi}}(y)\right\rangle\right\vert&=\left\vert \left\langle \pi_{\vec{n}}\left( \frac{\partial \vec{H}}{\partial\rho }\right) ,  \int_{\Gamma_{y,z}} d{\vec{\Phi}} \right\rangle \right\vert \\
&=\left\vert\int_{\Gamma_{y,z}}  \left\langle \pi_{\vec{n}(y)}\left( \frac{\partial \vec{H}}{\partial\rho } (y) \right) -\pi_{\vec{n}(x)}\left( \frac{\partial \vec{H}}{\partial\rho } (y) \right) ,  d{\vec{\Phi}}(x) \right\rangle \right\vert\\
&\leq \vert z\vert^2 \Vert \nabla \pi_{\vec{n}}\Vert_{L^\infty(\partial B(0,\vert z\vert))}  \left\Vert \frac{\partial \vec{H}}{\partial\rho }\right\Vert_{L^\infty(\partial B(0,\vert z\vert))}  e^{\lambda_{\vert z\vert}} .
\end{split}
\ee
Thanks to (\ref{pi}), we have, for all $\vec{w} \in \R^m$, that

$$\nabla \pi_{\vec{n}}(\vec{w})=(-1)^{m-1} (\nabla \vec{n})\res(\vec{n}\res \vec{w}) +(-1)^{m-1} \vec{n}\res((\nabla \vec{n}) \res \vec{w}).$$
Hence we have 
\beq
\label{c02}
\left\Vert\left\langle \pi_{\vec{n}}\left( \frac{\partial \vec{H}}{\partial\rho }\right) (y) , {\vec{\Phi}}(z)-{\vec{\Phi}}(y)\right\rangle\right\Vert_{L^\infty(\partial B(0,\vert z\vert))}\leq \vert z\vert^2 \Vert \nabla \vec{n}\Vert_{L^\infty(\partial B(0,\vert z\vert))}  \left\Vert \nabla  \vec{H} \right\Vert_{L^\infty(\partial B(0,\vert z\vert))}  e^{\lambda_{\vert z\vert}} .
\eeq
Combining  (\ref{c00}), (\ref{pit}), (\ref{c01}) and (\ref{c02}) we get

\be
\begin{split}
 \vert 2\pi C_0(z) \vert &\leq C \left( \vert z\vert \Vert H \Vert_{L^\infty(\partial B(0,\vert z\vert))}  \Vert \nabla \vec{n} \Vert_{L^\infty(\partial B(0,\vert z\vert))}    \Vert {\vec{\Phi}} -{\vec{\Phi}}(z) \Vert_{L^\infty(\partial B(0,\vert z\vert))}  \right. \\
 &+\left. \vert z\vert^3 \Vert \nabla \vec{n}\Vert_{L^\infty(\partial B(0,\vert z\vert))}  \left\Vert \nabla  \vec{H} \right\Vert_{L^\infty(\partial B(0,\vert z\vert))}  e^{\lambda_{\vert z\vert}}\right).  
 \end{split}
 \ee
Then using (\ref{er1}) and (\ref{er2}), we get
\be
 \vert 2\pi C_0(z) \vert  \leq C \Vert \nabla  \vec{n} \Vert_{L^{2}(B(0,2\vert z\vert) \setminus B(0,\frac{\vert z\vert}{2}))}^2 .
 \ee
Hence, we finally get 
$$ \vert 2\pi C_0(z) \vert \leq C \left( \int_{B(0,2\vert z\vert )\setminus B(0,\frac{\vert z\vert}{2})} \vert \nabla \vec{n}\vert^2 \, dz \right) \hbox{ for all } z \in B(0,1/2)\setminus B(0, 2r)
$$
Then, applying the following lemma, we get the desired result.\hfill$\square$

\begin{lemma}
\label{l21l}
Let $s:B(0,1/2)\setminus B(0,2r) \rightarrow \R_+$ such that 
$$s(z)=\frac{1}{\vert z\vert }  \left( \int_{B(0,2\vert z\vert )\setminus B(0,\vert z\vert /2) } \vert \nabla \vec{n} \vert^2 \; dz\right)^\frac{1}{2} $$
 Then $s^2(z) \vert z\vert $ is uniformly bounded in $L^{2,1}$ independently of $r$. \hfill $\Box$
\end{lemma}

\noindent{\bf Proof of lemma~\ref{l21l}.}

$$\frac{\partial (s^2(z) \vert z\vert ) }{\partial \rho} = -\frac{1}{\vert z\vert^2} \int_{B(0,2\vert z\vert )\setminus B(0,\vert z\vert /2) } \vert \nabla \vec{n} \vert^2 \; dz +\frac{1}{\vert z\vert }\int_{\partial B(0,2\vert z\vert )\cup \partial B(0,\vert z\vert /2) } \vert \nabla \vec{n} \vert^2 \; dz$$
In order to control the second term we can apply (\ref{er1}), which gives

$$ \frac{1}{\vert z\vert }\int_{\partial B(0,2\vert z\vert )\cup \partial B(0,\vert z\vert /2) } \vert \nabla \vec{n} \vert^2 \; dz \leq C\frac{1}{\vert z \vert^2 } \int_{B(0,4\vert z\vert )\setminus B(0,\vert z\vert /4) } \vert \nabla \vec{n} \vert^2 \; dz$$
Hence
$$ \int_{4r}^{1/4} \left\vert \frac{\partial (s^2(z) \vert z\vert ) }{\partial \vert z\vert} \right\vert \vert z\vert d\vert z\vert \leq  C\int_{B(0,1 )\setminus B(0,r) } \vert \nabla \vec{n} \vert^2 \; dz,$$
which gives the desired estimate thanks to the standard embedding of $W^{1,1}$ into $L^{2,1}$, see theorem 3.3.10 of \cite{Helein}.\hfill$\square$
\subsubsection{Estimating the second residue $\vec{C}_1$.}
\label{c1estsec}
\begin{lemma}
\label{c1est}
There exists $C>0$, independent of $r$ and $k$, such that
\beq
\label{estc1}
\left\Vert \frac{\vec{C}_1}{\rho} \right\Vert_{L^{2}(B(0,\frac{1}{2})\setminus B(0,2r)} \leq C.
\eeq
\hfill $\Box$
\end{lemma}

\noindent{\bf Proof of lemma~\ref{c1est}.}\\

Integrating by parts  (\ref{c1}) and using Stokes's theorem, we get for every $\rho\in (\eps, 1)$

\beq
\begin{split}
2\pi \vec{C}_1(z) &=\int_{\partial B(0,\rho)}  \partial_\tau \vec{L}\wedge {\vec{\Phi}} -(-1)^{m-1}2 (\star(\vec{n}\res \vec{H}))\res  \partial_\tau {\vec{\Phi}} \, d\sigma-\frac{1}{\rho} \int_{\partial B(0,\rho)} \vec{c}\wedge{\vec{\Phi}}  \, d\sigma +2\pi \vec{c}\wedge {\vec{\Phi}} \\
&=\int_{\partial B(0,\rho)}  - \vec{T}_\nu \wedge {\vec{\Phi}} +\frac{\vec{c}\wedge {\vec{\Phi}}}{\rho} -(-1)^{m-1}2 (\star(\vec{n}\res \vec{H}))\res  \partial_\tau {\vec{\Phi}} \, d\sigma-\frac{1}{\rho} \int_{\partial B(0,\rho)} \vec{c}\wedge{\vec{\Phi}}  \, d\sigma +2\pi \vec{c}\wedge {\vec{\Phi}} \\
&=\int_{\partial B(0,\rho)}  \vec{T}_\nu \wedge ({\vec{\Phi}}(z)-{\vec{\Phi}}(y))  -(-1)^{m-1} 2(\star(\vec{n}\res \vec{H}))\res  \partial_\tau {\vec{\Phi}}\, d\sigma ,
\end{split}
\eeq

where, as in the previous section,
$$\vec{T}_{\nu} =    \frac{\partial \vec{H}}{\partial\rho } -3\pi_{\vec{n}}\left( \frac{\partial \vec{H}}{\partial\rho }\right) - \star\left(\frac{\partial \vec{n}}{\partial \tau}\wedge \vec{H}\right)=  \pi_{T}\left( \frac{\partial \vec{H}}{\partial\rho } \right)-2\pi_{\vec{n}}\left( \frac{\partial \vec{H}}{\partial\rho }\right) - \star\left(\frac{\partial \vec{n}}{\partial \tau}\wedge \vec{H}\right). $$ 
Then 

$$2\pi \vec{C}_1(z) =\int_{\partial B(0,\rho)} \left( \pi_{T}\left( \frac{\partial \vec{H}}{\partial\rho } \right)-2\pi_{\vec{n}}\left( \frac{\partial \vec{H}}{\partial\rho }\right) - \star\left(\frac{\partial \vec{n}}{\partial \tau}\wedge \vec{H}\right)\right)\wedge ({\vec{\Phi}}(z)-{\vec{\Phi}}(y)) -(-1)^{m-1}2(\star(\vec{n}\res \vec{H}))\res   \partial_\tau {\vec{\Phi}}   \, d\sigma .$$ 
 Then, we easily check that except the second  term and the last term, i.e.

$$I=-2 \int_{\partial  B(0,\rho)}  \pi_{\vec{n}}\left( \frac{\partial \vec{H}}{\partial\rho }\right) \wedge({\vec{\Phi}}(z)-{\vec{\Phi}}(y)) \, d\sigma $$
and

$$II=-(-1)^{m-1}2 \int_{\partial  B(0,\rho)}   (\star(\vec{n}\res \vec{H}))\res \partial_\tau {\vec{\Phi}}  \, d\sigma ,$$
all the other terms can be controlled by $2\pi \rho e^{\lambda_\rho}\Vert \nabla \vec{n} \Vert_{L^\infty(\partial B(0,\rho))} \Vert H\Vert_{L^\infty(\partial B(0,\rho))} $ which is controlled by $\left( \int_{B(0,2\rho)\setminus B(0,\frac{\rho}{2})} \vert \nabla \vec{n}\vert^2 \, dz \right)$, as already done in the previous section for $C_0$.\\

For $I$ and $II$, we won't be able to derive a similar estimate than the one for $C_0$. The best we can do, applying  (\ref{er1}) and (\ref{er2}), is

$$I+II \leq  C \left( \int_{B(0,2\rho)\setminus B(0,\frac{\rho}{2})} \vert \nabla \vec{n}\vert^2 \, dz \right)^\frac{1}{2} .$$
Hence we get 

\beq
\label{finalc1}
\vert \vec{C}_1(z) \vert \leq  C \left( \int_{B(0,2\vert z \vert)\setminus B(0,\frac{\vert z \vert}{2})} \vert \nabla \vec{n}\vert^2 \, dz \right)^\frac{1}{2} \hbox{ for all }z \in B(0,1/2)\setminus B(0,2r) ,
\eeq
which implies that the desired estimate.\hfill$\square$\\

\section{Proof of theorem~\ref{th-main}.}

After deriving some equations relating $\vec{H}$ together with $S$ and $\vec{R}$, we prove that an appropriate combination of those quantities is bounded in $L^{2,1}$.
This will be done using some additional conservation laws combined with results from integrability by compensation theory from \cite{LR} or proved in the appendix. Finally combining those estimates with some algebraic fact, we prove that $\vec{H}$ converges to $0$ in $L^{2}$ norm in the neck region.
\subsection{Equations satisfied by $\vec{H}$, $\vec{R}$ and $S$.}
Equation  (\ref{e2}) gives
$$ \nabla \vec{R}  \res \nabla^\bot {\vec{\Phi}}= (\vec{L}\wedge \nabla {\vec{\Phi}})\res \nabla^\bot {\vec{\Phi}} +(-1)^{m-1}2(\star(\vec{n}\res \vec{H})\res \nabla {\vec{\Phi}})\res \nabla^\bot {\vec{\Phi}}+ \vec{C}_1 \nabla^\bot \log(\rho) \res \nabla^\bot {\vec{\Phi}}  ,$$
which simplify into
$$ \nabla \vec{R}  \res \nabla^\bot {\vec{\Phi}}= \langle \vec{L}, \nabla^\bot {\vec{\Phi}}\rangle \nabla {\vec{\Phi}} +(-1)^{m-1}2(\star(\vec{n}\res \vec{H})\res \nabla {\vec{\Phi}})\res \nabla^\bot {\vec{\Phi}}+ \frac{1}{\rho} \vec{C}_1\res \frac{\partial {\vec{\Phi}}}{\partial\rho } .$$
Since, see X.210 of \cite{conf}, for any normal vector $\vec{N}$, 
\beq
\label{X.210}
\star(\vec{n}\res \vec{N})=(-1)^{m-1}\vec{e}_1 \wedge \vec{e}_2\wedge \vec{N},
\eeq
we deduce that 
$$ (-1)^{m-1}2(\star(\vec{n}\res \vec{H})\res \nabla {\vec{\Phi}})\res \nabla^\bot {\vec{\Phi}} = -\,4\ e^{2\lambda} \ \vec{H}.$$
Moreover (\ref{e1}) gives
$$ \langle \vec{L}, \nabla^\bot {\vec{\Phi}}\rangle =\nabla^\bot S +C_0 \nabla \log(\rho),$$
which finally implies
\beq
\label{Heq}
4\ e^{2 \lambda}\ \vec{H}=- \nabla \vec{R} \res \nabla^\bot {\vec{\Phi}} + \nabla^\bot S \cdot \nabla {\vec{\Phi}}  + \frac{1}{\rho} C_0 \frac{\partial {\vec{\Phi}}}{\partial\rho }   + \frac{1}{\rho} \vec{C}_1\res \frac{\partial {\vec{\Phi}}}{\partial\rho }.
\eeq
We now compute an autonomous system satisfies by $\vec{R}$ and $S$. From (\ref{e1}) and (\ref{e2}), we have
\beq
\left\{\begin{array}{l} \nabla^\bot S = \langle \vec{L}, \nabla^\bot {\vec{\Phi}}\rangle -C_0\nabla  \log(\rho) \\[5mm]
\nabla^\bot \vec{R} =  \vec{L}\wedge \nabla^\bot {\vec{\Phi}} + (-1)^{m-1}2(\star(\vec{n}\res \vec{H})\res \nabla^\bot {\vec{\Phi}})-\vec{C}_1 \nabla \log(\rho)\end{array}\right.
\eeq
But we get, for any normal vector $\vec{N}$,  as a consequence of (\ref{X.210}), 
\beq
\label{X.211}
(-1)^{m-1}(\star(\vec{n}\res \vec{N}))\res \nabla {\vec{\Phi}}=-\nabla^\bot {\vec{\Phi}} \wedge \vec{N}. 
\eeq
Then, combining  (\ref{e2}),  (\ref{X.211}) and the definition of $\bullet$, see (\ref{b1}) and (\ref{b2}), we deduce that
\beq 
\label{XX}
\vec{n}\bullet \nabla^\bot \vec{R} =(\vec{n}\res \pi_{\vec{n}} (\vec{L}))\wedge \nabla^\bot {\vec{\Phi}} - 2(\vec{n}\res \vec{H})\wedge \nabla {\vec{\Phi}} -\vec{n}\bullet \vec{C}_1 \nabla \log(\rho)  .
\eeq
From X.215 of \cite{conf} , if $\vec{N}$ is a normal vector field, we get

\beq 
\star( (\vec{n}\res \vec{N} )\wedge \nabla^\bot {\vec{\Phi}})=(-1)^{m-1} \nabla {\vec{\Phi}} \wedge \vec{N}
\eeq
which gives, combining it with (\ref{XX}),

\beq 
\label{XXX}
\star(\vec{n}\bullet \nabla^\bot  \vec{R}) =(-1)^m \pi_{\vec{n}} (\vec{L})\wedge \nabla {\vec{\Phi}}  + (-1)^{m-1}2 \nabla^\bot {\vec{\Phi}} \wedge \vec{H}   - \star(\vec{n}\bullet \vec{C}_1 \nabla  \log(\rho)) .
\eeq
Moreover, combining (\ref{e2}) and (\ref{X.211}), we get 
\be
\begin{split}
\nabla \vec{R} &=  \vec{L}\wedge \nabla {\vec{\Phi}} +2(\star(\vec{n}\res \vec{H})\res \nabla {\vec{\Phi}})+\vec{C}_1 \nabla^\bot \log(\rho) \\[3mm]
&= \pi_{\vec{n}}(\vec{L})\wedge \nabla {\vec{\Phi}} -2 \nabla^\bot {\vec{\Phi}} \wedge \vec{H}+ \pi_{T}(\vec{L})\wedge \nabla {\vec{\Phi}} +\vec{C}_1 \nabla^\bot \log(\rho)
\end{split}
\ee
Then thanks to (\ref{XXX}), we get
\beq 
\label{X4}
(-1)^{m} \star(\vec{n}\bullet \nabla^\bot  \vec{R}) =\nabla \vec{R} - \pi_{T} (\vec{L})\wedge \nabla {\vec{\Phi}}   -\vec{C}_1 \nabla^\bot \log(\rho)  - (-1)^m \star(\vec{n}\bullet \vec{C}_1 \nabla  \log(\rho)) .
\eeq
We easily check, using (\ref{e1}), that
$$\pi_{T} (\vec{L})\wedge \nabla {\vec{\Phi}}= \nabla^\bot S\star \vec{n} + C_0 \nabla \log(\rho) \star \vec{n} ,$$
which finally gives with (\ref{X4})
\beq
\label{main1} 
\begin{split}
\nabla \vec{R}&= (-1)^{m} \star(\vec{n}\bullet \nabla^\bot  \vec{R}) + \nabla^\bot  S \star \vec{n} + C_0 \nabla \log(\rho) \star \vec{n} \\[3mm] 
& +\vec{C}_1\, \nabla^\bot \log(\rho) +(-1)^m \star(\vec{n}\bullet \vec{C}_1 \nabla  \log(\rho)).
\end{split}
\eeq
Now, taking the scalar product again with $\star \vec{n}$ we get 
\beq
\label{main11}
\nabla S= -\langle \nabla^\bot \vec{R},\star \vec{n} \rangle   + C_0 \nabla^\bot \log(\rho) - \langle \vec{C}_1\, \nabla \log(\rho),\star\vec{n}\rangle ,
\eeq
here we used the fact that 
$$\langle \star(\vec{n}\bullet \nabla^\bot  \vec{R}) , \star \vec{n} \rangle =- \langle  \star(\vec{n}\bullet \vec{C}_1 \nabla  \log(\rho)), \star \vec{n}\rangle   $$
since, in the right hand side of  (\ref{XXX}), excepted the term involving $\vec{C}_1$, all terms are a  linear combination of wedges of tangent and normal vectors.\\

Taking divergence we finally get 
\beq
\label{main2}
\left\{\begin{array}{l}

\Delta \vec{R}= (-1)^{m} \star(\nabla \vec{n}\bullet \nabla^\bot  \vec{R}) + \nabla^\bot  S \nabla \star \vec{n} +  \mbox{div}(C_0 \nabla \log(\rho)  \star \vec{n} ) \\[3mm]
\quad\quad+ \mbox{div}(\vec{C}_1 \nabla^\bot \log(\rho)) +(-1)^m \mbox{div}(\star(\vec{n}\bullet \vec{C}_1 \nabla \log(\rho)))

\\[5mm]
\Delta  S= -\langle \nabla^\bot \vec{R},\nabla \star \vec{n} \rangle + \mbox{div}( C_0 \nabla^\bot \log(\rho)) \\[3mm]
\quad\quad- \mbox{div}(\langle \vec{C}_1 \nabla \log(\rho),\star\vec{n}\rangle)
\end{array}\right.
\eeq

\subsection{Estimate I of $\vec{C}_1$}
\begin{lemma}
\label{c1e1}
There exists $C>0$, independent of $r$, such that
\beq
\label{est1c1}
\Vert \langle \vec{C}_1 \nabla \log(\rho),\star\vec{n}\rangle \Vert_{L^{2,1}(B(0,\frac{1}{2})\setminus B(0,2r)} \leq C.
\eeq
\hfill $\Box$
\end{lemma}

\noindent{\bf Proof of lemma~\ref{c1e1}.}
Proceeding as in  lemma \ref{c1est}, we easily  control every term except the one involving $I$ and $II$. That is to say, for all $z\in \partial B(0,\rho)$ with $2r<\rho<1/2$, we get 

\beq
\label{estc11}
\begin{split}
\vert \langle \vec{C}_1(z),\star\vec{n}(z)\rangle\vert &=O \left( \left\vert  \int_{\partial  B(0,\rho)} \left\langle \star \vec{n}(z), \left(    2\pi_{\vec{n}} \left(\frac{\partial \vec{H}}{\partial\rho }\right)  \wedge ({\vec{\Phi}}(z)-{\vec{\Phi}}(y))  +2  (\star(\vec{n}\res \vec{H}))\res  \partial_\tau {\vec{\Phi}}  \right)\right\rangle  \, d\sigma_z\right\vert\right)  \\
 &+ O  \left( \int_{B(0,2\rho)\setminus B(0,\frac{\rho}{2})} \vert \nabla \vec{n}\vert^2 \, dz  \right) \\
  &=O \left( \left\vert  \int_{\partial  B(0,\rho)} \left\langle \star (\vec{n}(z) -\vec{n}(y)), \left(  2\pi_{\vec{n}} \left(\frac{\partial \vec{H}}{\partial\rho }\right)  \wedge ({\vec{\Phi}}(z)-{\vec{\Phi}}(y))  +2(-1)^m \frac{\partial {\vec{\Phi}}}{\partial \rho}  \wedge \vec{H} \right)\right\rangle  \, d\sigma_z\right\vert\right)  \\
 &+ O \left( \left\vert  \int_{\partial  B(0,\rho)} \left\langle \star \vec{n}(y), \left(  2\pi_{\vec{n}} \left(\frac{\partial \vec{H}}{\partial\rho }\right)  \wedge({\vec{\Phi}}(z)-{\vec{\Phi}}(p))   +2(-1)^m \frac{\partial {\vec{\Phi}}}{\partial \rho}  \wedge \vec{H} \right)\right\rangle  \, d\sigma_z\right\vert\right)  \\
 &+O \left( \int_{B(0,2\rho)\setminus B(0,\frac{\rho}{2})} \vert \nabla \vec{n}\vert^2 \, dz  \right)  
 \end{split}
\eeq

Here we used (\ref{X.211}) to simplify the last term of the first line. Proceeding as in the estimate of $C_0$, we get
 
\beq
\label{estc11suite}
\begin{split}
\vert \langle \vec{C}_1(z),\star\vec{n}(z)\rangle\vert &= O \left( \vert z\vert^2 \Vert \nabla \vec{n} \Vert_{L^\infty(\partial B(0,\rho))} \left[\vert z\vert \Vert \nabla \vec{H} \Vert_{L^\infty(\partial B(0,\rho))} +\Vert \vec{H}\Vert_{L^\infty(\partial B(0,\rho))} \right]e^{\lambda_\rho}\right)  \\
&+O \left( \left\vert  \int_{\partial  B(0,\rho)} \left\langle \star \vec{n}(y), \left(2 \pi_{\vec{n}} \left(\frac{\partial \vec{H}}{\partial\rho }\right)  \wedge ({\vec{\Phi}}(z)-{\vec{\Phi}}(y))  +2(-1)^m \frac{\partial {\vec{\Phi}}}{\partial \rho}  \wedge \vec{H} \right)\right\rangle  \, d\sigma_z\right\vert\right)  \\
&+O \left( \int_{B(0,2\rho)\setminus B(0,\frac{\rho}{2})} \vert \nabla \vec{n}\vert^2 \, dz  \right)  \\
 &=O\left( \int_{B(0,2\rho)\setminus B(0,\frac{\rho}{2})} \vert \nabla \vec{n}\vert^2 \, dz  \right) .
\end{split}
\eeq
Here we use the fact that 
$$\left\langle \star \vec{n}(y), \left( 2\pi_{\vec{n}} \left(\frac{\partial \vec{H}}{\partial\rho }\right)  \wedge ({\vec{\Phi}}(z)-{\vec{\Phi}}(y))  +2(-1)^m \frac{\partial {\vec{\Phi}}}{\partial \rho}  \wedge \vec{H} \right)\right\rangle=0,$$
since $\star \vec{n}=\vec{e}_{1}\wedge \vec{e}_2$ and that  $\pi_{\vec{n}} \left(\frac{\partial \vec{H}}{\partial\rho }\right)  \wedge ({\vec{\Phi}}(z)-{\vec{\Phi}}(y))  +2(-1)^m \frac{\partial {\vec{\Phi}}}{\partial \rho}  \wedge \vec{H}$ can easily be rewritten as a linear combination of $\vec{n}^\alpha\wedge \vec{n}^\beta$ and  $\vec{n}^\alpha\wedge \vec{e}_i$. Finally, using (\ref{estc11suite}),  we conclude as in the proof of lemma \ref{c0est}.\hfill$\square$

\subsection{Estimate II of $\vec{C}_1$.}

\begin{lemma}
\label{c1e3}
There exists $C>0$, independent of $r$, such that
\beq
\label{est2c1}
\Vert \nabla \vec{C}_1 \Vert_{L^{2}(B(0,\frac{1}{2})\setminus B(0,2r)} \leq C.
\eeq
\hfill $\Box$
\end{lemma}
\noindent{\bf Proof of lemma~\ref{c1e3}.} First we observe that $\nabla \vec{C}_1 =\vec{c}\wedge \nabla {\vec{\Phi}}$, hence using the definition of $\vec{c}$, we get 

\beq
\label{gradc1}
\begin{split}
\vert \nabla \vec{C}_1 (z) \vert &\leq C \vert z\vert e^{\lambda_{\vert z\vert}}(  \vert \nabla H\vert  + e^{-\lambda_{\vert z\vert}} \vert \nabla \vec{n}\vert^2 )\\
&\leq C \frac{1}{\vert z\vert} \left( \int_{B(0,2\vert z\vert)\setminus B(0,\frac{\vert z\vert}{2})} \vert \nabla \vec{n}\vert^2 \, dz  \right) ^\frac{1}{2}
\end{split}
\eeq
The first inequality is a  consequence of (\ref{defc})  and the Harnack estimate on the conformal factor, while the second uses  (\ref{er1}), (\ref{er2}) and the fact that  the norm on the dyadic annuli  are uniformly bounded. As already observe the right hand side is uniformly bounded into $L^2$.\hfill$\square$

\begin{lemma}
\label{c1e4}
There exists $C>0$, independent of $r$, such that
\beq
\label{est2c1}
\vert \Delta \vec{C}_1(z) \vert \leq \frac{C}{\vert z\vert^2}   \left( \int_{B(0,2\vert z\vert)\setminus B(0,\frac{\vert z\vert}{2})} \vert \nabla \vec{n}\vert^2 \, dz\right)  \hbox{ for all } z\in B(0,1/2)\setminus B(0,2r),
\eeq
In particular, $ \Delta \vec{C}_1$ is uniformly bounded in $L^1$.
\hfill $\Box$
\end{lemma}
\noindent{\bf Proof of lemma~\ref{est2c1}.}
We have $\Delta \vec{C}_1= \vec{c}\wedge \Delta {\vec{\Phi}} =2e^{2\lambda}\vec{c}\wedge \vec{H}$, which gives, for all $z\in B(0,1/2)\setminus B(0,2r)$,

$$\vert \Delta \vec{C}_1(z) \vert \leq C \left\vert e^{2\lambda} \vec{H}(z) \wedge \int_{\partial B(0,\vert z\vert)}  \left(\partial_\nu \vec{H}-3\pi_{\vec{n}}(\partial_\nu \vec{H}) - \star(\partial_\tau \vec{n}\wedge \vec{H}\right)\, d\sigma \right\vert$$
Remembering (\ref{pit}), we have
$$\vert \Delta \vec{C}_1(z) \vert \leq C \left\vert e^{2\lambda} \vec{H}(z) \right\vert_{L^\infty(\partial B(0,\vert z\vert))}  \vert z\vert   \left( \Vert \nabla \vec{H} \Vert_{L^\infty(\partial B(0,\vert z\vert))}  +\Vert \nabla \vec{n}\Vert_{L^\infty(\partial B(0,\vert z\vert))} \Vert \vec{H}\Vert_{L^\infty(\partial B(0,\vert z\vert))}  \right).$$
Combining this with (\ref{er1}) and (\ref{er2}), we get
$$\vert \Delta \vec{C}_1(z) \vert \leq \frac{C}{\vert z\vert^2}   \left( \int_{B(0,2\vert z\vert)\setminus B(0,\frac{\vert z\vert}{2})} \vert \nabla \vec{n}\vert^2 \, dz\right)  .$$
\hfill$\square$
\begin{lemma}
\label{c1e5}
There exists $C>0$, independent of $r$, such that
\beq
\label{est2c1}
\Vert \nabla \vec{C}_1 \Vert_{L^{2}(B(0,\frac{1}{2})\setminus B(0,2r))}  \leq  C \sup_{\rho\in (r,1)} \left( \int_{B(0,2\rho)\setminus B(0,\frac{\rho}{2})} \vert \nabla \vec{n}\vert^2 \, dz\right)^\frac{1}{4}
\eeq
\hfill $\Box$
\end{lemma}
\noindent{\bf Proof of lemma~\ref{c1e5}.}
A simple integration by part gives
$$ \Vert \nabla \vec{C}_1 \Vert_2^2 \leq  \left\vert  \int_{B(0,\frac{1}{2})\setminus B(0,2r))}  \langle \vec{C}_1 \Delta \vec{C}_1 \, dz +  \int_{\partial B(0,\frac{1}{2})\cup \partial B(0,2r))}  \langle \vec{C}_1 \partial_{\nu} \vec{C}_1 \rangle \, d\sigma\right\vert $$
Now using (\ref{finalc1}), (\ref{gradc1}) and lemma~\ref{c1e4} we obtain the desired result.\hfill$\square$

\subsection{ $L^{2,\infty}-$bounds on  $\nabla \vec{R}$  and $\nabla S$.}
Using (\ref{e1}) and (\ref{e2})  we get

\beq
\label{L2iS}
\vert \nabla S\vert = O \left( \vert \vec{L}\vert e^{\lambda_\rho} +\left\vert \frac{C_0}{\rho}\right\vert \right)  \hbox{ on } \D\setminus B(0,r)
\eeq
and
\beq
\label{L2iR}
\vert \nabla \vec{R} \vert = O \left( \vert \vec{L}\vert e^{\lambda_\rho} + \vert \vec{H} \vert e^{\lambda_\rho} + \left\vert \frac{\vec{C}_1}{\rho}\right\vert \right) \hbox{ on } \D\setminus B(0,r).
\eeq
Now, combining (\ref{L2iS}), (\ref{L2iR}), (\ref{eL}), (\ref{er2}), lemma \ref{c0est} , lemma \ref{c1est}  and (\ref{r1}), we get
\beq
\label{L2i}
\Vert \nabla S\Vert_{L^{2,\infty}(\D \setminus B(0,r)} +\Vert \nabla \vec{R} \Vert_{L^{2,\infty}(\D \setminus B(0,r))} \leq C,
\eeq
with $C$ independent of $r$.
\subsection{ $L^2-$bounds on  $\nabla S$  and $\nabla \vec{R}$ .}
\label{L2}
Let denote $S_\rho=\frac{1}{\vert \partial B(0,\rho)\vert} \int_{\partial B(0,\rho)} S\; d\sigma$ and $\vec{R}_\rho=\frac{1}{\vert \partial B(0,\rho)\vert} \int_{\partial B(0,\rho)} \vec{R}\; d\sigma$. Thanks to equation (\ref{main2}) and the following  lemma it suffies to control $\int_r^1 \left\vert \frac{d S_\rho }{d\rho}  \right\vert^2 \; \rho d\rho $  and  $\int_r^1 \left\vert \frac{d \vec{R}_\rho}{d\rho}   \right\vert^2 \; \rho d\rho $ to  prove that  $\nabla S$  and $\nabla \vec{R}$ are uniformly bounded in $L^2(B(0,1/2)\setminus B(0,2r))$.
\begin{lemma}[lemma 2.4 \cite{LR}]
\label{LR1}
Let  $a,b \in L^2(B_1)$,  $0<\eps<\frac{1}{4}$,  assume that $\nabla a\in L^{2,\infty}(B_1)$ and that $\nabla b\in L^2(B_1)$,  let $\varphi\in W^{1,(2,\infty)}(B_1\setminus B_\eps)$ a solution of 
\beq
\Delta\varphi = a_xb_y-a_yb_x \hbox{ on } B_1\setminus B_\eps ,\\
\eeq
Denote, for $\eps\le r\le1$, $\varphi_\rho:=(2\pi\, \rho)^{-1}\,\int_{\p B(0,\rho)}{\varphi}\, d\sigma$ and assume
\beq
\label{x1}
\int_{\eps}^1\left\vert \frac{d \varphi_\rho}{d\rho} \right\vert^2\ r\ dr<+\infty\quad.
\eeq
Then, for $0\leq \delta<1$, there exists a positive constant $C(\delta)>0$ independent of $\eps$ and $\varphi$ such that
\beq
\label{x2}
\begin{split}
\ds\Vert \nabla \varphi \Vert_{L^{2}\left(B_\delta\setminus  B_{\delta^{-1}\eps }\right)} &\leq C(\delta) \left( \ \Vert \nabla a\Vert_{{2,\infty}}\,  \Vert \nabla b\Vert_{2} + \|\nabla\varphi_\rho\|_{L^2(B_1\setminus B_\eps)} \right. \\
& \left.+ \ \Vert\nabla \varphi\Vert_{L^{2,\infty}(B_1\setminus B_\eps)}\right).
\end{split}
\eeq
\hfill $\Box$
\end{lemma}
Thanks to (\ref{main11}) and (\ref{main1}) we get 

\beq
\label{eS1}
\frac{d S_\rho}{d\rho}    = \int_0^{2\pi} \left\langle \frac{1}{\rho} \frac{\partial \vec{R}}{ \partial \theta },\star \vec{n} \right\rangle\; d\theta  -  \int_0^{2\pi}   \left\langle \star \vec{n}, \frac{\vec{C}_1}{\rho} \right\rangle\; d\theta  
\eeq
and
\beq
\frac{d \vec{R}_\rho}{d\rho}   =  - \int_0^{2\pi} (-1)^{m} \star\left(\vec{n}\bullet \frac{1}{\rho} \frac{\partial \vec{R}}{\partial \theta}\right) \; d\theta -  \int_0^{2\pi}  \frac{1}{ \rho} \frac{\partial  S}{\partial \theta} \star \vec{n} \; d\theta + \int_0^{2\pi} \frac{C_0 \star \vec{n}}{\rho}   \; d\theta  + \int_0^{2\pi} (-1)^m \star(\vec{n}\bullet \frac{\vec{C}_1 }{\rho}) d\theta .
\eeq
On the one hand 
\be
\left\vert  \int_0^{2\pi} \left\langle \frac{1}{\rho} \frac{\partial \vec{R}}{ \partial \theta },\star \vec{n} \right\rangle\; d\theta  \right\vert =  \left\vert \int_0^{2\pi} \left\langle\frac{1}{\rho}  \frac{\partial \vec{R}}{ \partial \theta },\star (\vec{n}(z)-\vec{n}(\vert z\vert,0))\right \rangle\; d\theta  \right\vert \leq \rho \Vert \nabla \vec{n} \Vert_{L^{\infty} (\partial B(0,\rho))}  \Vert \nabla \vec{R} \Vert_{L^{\infty} (\partial B(0,\rho))}  .
\ee
On the other hand, thanks to (\ref{e2}),  
\beq
\label{eS2}
\Vert \nabla \vec{R} \Vert_{L^{\infty} (\partial B(0,\rho))}  \leq e^{\lambda_\rho} \left(\Vert  \vec{L} \Vert_{L^{\infty} (\partial B(0,\rho))}  + \Vert  \vec{H} \Vert_{L^{\infty} (\partial B(0,\rho))}\right) +\left\Vert \frac{\vec{C}_1}{\rho}\right\Vert_{L^{\infty} (\partial B(0,\rho))}.
\eeq
Then combining (\ref{er1}), (\ref{eL}), (\ref{finalc1}), (\ref{eS1}), (\ref{eS2}) we finally get 
$$ \left\vert \frac{d S_\rho}{d\rho}   \right\vert^2 \ \leq \frac{C}{\rho^2} \int_{B(0,2\rho)\setminus B(0,\rho/2)} \vert \nabla \vec{n} \vert^2  \; dz  $$
which gives
$$\int_r^1 \left\vert \frac{d S_\rho}{d\rho}   \right\vert^2 \; \rho d\rho \leq C \int_{\D\setminus B(0,r)} \vert \nabla \vec{n} \vert^2  \; dz  .$$
Mutatis mutandis, we get
$$\int_r^{1} \left\vert \frac{d  \vec{R}_\rho}{d\rho}  \right\vert^2\;\rho  d\rho \leq C \int_{\D\setminus B(0,r)} \vert \nabla \vec{n} \vert^2  \; dz  .$$
Now applying lemma \ref{LR1}, we get that
\be
\Vert \nabla   (S-\psi_S) \Vert_{L^2(B(0,\frac{1}{2})\setminus B(0,2r))} \leq C \Vert \nabla \vec{n} \Vert_2 \Vert \nabla \vec{R}  \Vert_{2,\infty}  \ee
where $\psi_S \in W^{1,1}_0(\D\setminus B(0,r))$ satisfies
$$\Delta \psi_S=  div( C_0 \nabla^\bot \log(\rho))  - div(\langle \vec{C}_1 \nabla \log(\rho),\star\vec{n}\rangle) ,$$
and thanks to standard elliptic theory, we have  
$$\int_r^{1} \left\vert \frac{d}{d\rho} (\psi_{S})_\rho  \right\vert^2\rho \; d\rho \leq \Vert \nabla \psi_S\Vert_2^2 \leq  C  \left(  \left\Vert \frac{C_0}{\rho}\right\Vert_{L^2(\D\setminus B(0,r))}^2 + \left\Vert \frac{\vec{C}_1}{\rho}\right\Vert_{L^2(\D\setminus B(0,r))}^2\right)  $$
which is bounded thanks to  lemmas \ref{c0est} and \ref{c1est}.\\

We also get 
\be
\begin{split}
\Vert \nabla( \vec{R} -\psi_{\vec{R}})\Vert_{L^2(B(0,\frac{1}{2})\setminus B(0,2r))} &\leq C\left(\Vert \nabla \vec{n} \Vert_2 (\Vert \nabla \vec{R} \Vert_{L^{2,\infty}(\D\setminus B(0,r))} + \Vert \nabla S \Vert_{L^{2,\infty}(\D\setminus B(0,r))}\right) \\
&+\left\Vert \frac{C_0}{\rho}\right\Vert_{L^2(\D\setminus B(0,r))} + \left\Vert \frac{\vec{C}_1}{\rho}\right\Vert_{L^2(\D\setminus B(0,r))} 
\end{split}
 \ee
 where $\psi_{\vec{R}} \in W^{1,1}_0(\D\setminus B(0,r))$ satisfies
$$\Delta \psi_{\vec{R}}=    div(C_0 \nabla \log(\rho)  \star \vec{n} ) +  div(\vec{C}_1 \nabla^\bot \log(\rho)  )  +div((-1)^m \star(\vec{n}\bullet \vec{c}_1 \nabla \log(\rho))),$$
and thanks to standard elliptic theory, we have  
$$\int_r^{1} \left\vert \frac{d}{d\rho} (\psi_{\vec{R}})_\rho  \right\vert^2\rho \; d\rho \leq \Vert \nabla \psi_{\vec{R}}\Vert_2^2\leq  C  \left\Vert \frac{C_0}{\rho}\right\Vert_{L^2(\D\setminus B(0,r))}^2 + \left\Vert \frac{\vec{C}_1}{\rho}\right\Vert_{L^2(\D\setminus B(0,r))}^2 $$
which is bounded thanks to  lemma \ref{c0est} and lemma \ref{c1est}.\\  

Finally, using (\ref{L2i}), we can  conclude that 
 \beq
 \label{L2est}
\Vert \nabla   S\Vert_{L^2(B(0,\frac{1}{2})\setminus B(0,2r))}  +\Vert \nabla \vec{R} \Vert_{L^2(B(0,\frac{1}{2})\setminus B(0,2r))} \leq C, 
\eeq
with $C$ independent of $r$.
\subsection{$L^{2,1}-$bounds on $\nabla S$.}
First, integrating (\ref{main11}), we get
\beq
\label{deltan}
\begin{split}
\frac{d S_\rho}{d\rho}    &= \int_0^{2\pi} \langle \frac{\partial \vec{R}}{\rho \partial \theta },\star \vec{n}(z)\rangle\; d\theta  -  \int_0^{2\pi}   \langle \star \vec{n}, \vec{C}_1\frac{1}{\rho} \rangle\; d\theta   \\
&= \int_0^{2\pi} \langle \frac{\partial \vec{R}}{\rho \partial \theta },\star (\vec{n}(z)-\vec{n}((0,\rho))) \rangle\; d\theta  -  \int_0^{2\pi}   \langle \star \vec{n}, \vec{C}_1\frac{1}{\rho} \rangle\; d\theta   \\
&\leq \left(\int_0^{2\pi} \vert \vec{n}(z)-\vec{n}((0,\rho))\vert^2 \; d\theta \right)^\frac{1}{2} \left(\int_0^{2\pi} \vert \nabla \vec{R}\vert^2 \; d\theta \right)^\frac{1}{2}  + \left\vert \int_0^{2\pi} \langle \star \vec{n}, \vec{C}_1\frac{1}{\rho} \rangle\; d\theta \right\vert \\
&\leq \left(\rho \int_{\partial B(0,\rho)} \vert \nabla \vec{n}\vert^2 \; d\sigma \right)^\frac{1}{2} \left(\frac{1}{\rho}\int_{\partial B(0,\rho)}  \vert \nabla \vec{R}\vert^2 \; d\sigma \right)^\frac{1}{2}  \\
& +O\left( \frac{1}{\rho} \int_{B(0,2\rho)\setminus B(0,\frac{\rho}{2})} \vert \nabla \vec{n}\vert^2 \, dz  \right)  
\end{split}
\eeq
Here we used (\ref{estc11suite}) to estimate the last term of third line.  Finally, we get 
$$\int_r^1 \left\vert \frac{d S_\rho}{d\rho}   \right\vert \;  d\rho \leq C \int_{\D\setminus B(0,r)} \vert \nabla \vec{n} \vert^2 \; dz   \leq C $$
Since both $S$ and $\vec{R}$ are defined up to an additive constant, we can assume that $S_{r}=0$ and $\vec{R}_r=0$, which gives

$$ \vert S_\rho\vert \leq C \;\; \forall \rho\in (r,1) .$$
Now let us recall the following lemma from \cite{LR}.
\begin{lemma}
\label{lm-2.2}
Let  $a,b \in W^{1,2}(B_1)$,  $0<\eps<\frac{1}{4}$,  and $\varphi\in W^{1,1}(B_1\setminus B_\eps)$ a solution of 
\beq
\label{lbis}
\left\{
\begin{array}{l}
\displaystyle \Delta \varphi = a_xb_y-a_yb_x \hbox{ on } B_1\setminus B_\eps \\[5mm]
\displaystyle \int_{\partial B_\eps} \varphi\;d\sigma =0,\\[5mm]
\displaystyle \left\vert \int_{\partial B_1 } \varphi  \, d\sigma \right\vert \leq K,\\
\end{array}
\right.
\eeq
where $K$ is a constant independent of $\eps$. Then there exists a positive constant $C$ independent of $\eps$ such that 
$$\Vert \nabla \varphi\Vert_{L^{2,1}\left( B(0,1/2)\setminus   B(0,2\eps )\right)} \leq C( \Vert \nabla a\Vert_2\,  \Vert \nabla b\Vert_2 + \Vert \nabla \varphi \Vert_2 + 1)\quad .$$
\hfill $\Box$
\end{lemma}
Thanks to lemma \ref{c0est} and lemma \ref{c1e1}, we easily get that the solution $\psi_S \in W^{1,1}_0(\D\setminus B(0,r))$  of
$$\Delta \psi_S=   div(C_0 \nabla^\bot \log(\rho)\star\vec{n})-div(\langle \vec{C}_1 \nabla \log(\rho),\star\vec{n}\rangle) ,$$
satisfies the following estimate
\beq
\label{psiS}
\Vert\psi_S \Vert_{\infty} \leq  C\Vert \nabla \psi_S \Vert_{L^{2,1}(\D\setminus B(0,r))}  \leq C \left( \left\Vert C_0\nabla^\bot \log(\rho) \star \vec{n} \right\Vert_{L^{2,1}(\D\setminus B(0,r))} + \left\Vert \frac{\langle \vec{c}_1, \star \vec{n} \rangle }{\rho}\right\Vert_{L^{2,1}(\D\setminus B(0,r))} \right)  \leq C, 
\eeq
with $C$ independent of $r$. Here we have used the injection $W^{1,(2,1)}_0 \subset L^\infty$, see theorem 3.3.4 of \cite{Helein}.
Applying lemma \ref{lm-2.2}, we get 
\be
\Vert \nabla   (S-\psi_S)\Vert_{L^{2,1}(B(0,\frac{1}{2})\setminus B(0,2r))} \leq C(\Vert \nabla \vec{n} \Vert_2 (\Vert \nabla \vec{R} \Vert_{2} + \Vert \nabla S\Vert_2 +1))
\ee
Finally using again (\ref{psiS}), we get that $\Vert \nabla S\Vert_{L^{2,1}(B(0,\frac{1}{2})\setminus B(0,2r))} $ is uniformly bounded.

\subsection{Final estimate} 
Thanks to (\ref{Heq}), (\ref{finalc0}) and the previous paragraph, we get that 

\beq
\label{Hfin}
\left\Vert 4 e^{ \lambda} \vec{H}+ e^{-\lambda} \nabla \vec{R} \res \nabla^\bot{\vec{\Phi}} - e^{-\lambda} \frac{1}{\rho} \vec{C}_1\res \frac{\partial {\vec{\Phi}}}{\partial \rho } \right\Vert_{L^{2,1}(B(0,\frac{1}{2})\setminus B(0,2r))}  \leq C,
\eeq
with $C$ independent of $r$.\\

But (\ref{main1}) gives that 

$$\nabla \vec{R}-(-1)^m \star(\vec{n}\bullet \vec{C}_1 \nabla \log(\rho)) =  (-1)^{m} \star(\vec{n}\bullet \nabla^\bot  \vec{R}) + \nabla^\bot  S \star \vec{n} + C_0 \nabla \log(r) \star \vec{n} +\vec{C}_1 \nabla^\bot \log(\rho).$$
Let  $f,g \in W^{1,1}_0(\D\setminus B(0,r)$ such
\beq
\label{eqf}
\Delta \vec{f} =   \nabla \vec{C}_1 \nabla^\bot \log(\rho)
\eeq
and
$$ \Delta \vec{g} = -(-1)^m \star(\nabla^\bot (\vec{n}\bullet \vec{C}_1) \nabla \log(\rho)).$$
Thanks to lemma \ref{lm-a3}, lemma \ref{c1e5} and (\ref{finalc1}), we have
\beq
\label{fgest}
\Vert \nabla f\Vert_2 +\Vert \nabla g\Vert_2 \leq C  \sup_{\rho\in (r,1)} \left( \int_{B(0,2\rho)\setminus B(0,\frac{\rho}{2})} \vert \nabla \vec{n}\vert^2 \, dz\right)^\frac{1}{4} 
\eeq
with $C$ independent of $r$.\\

Then we set $\vec{X}=  \nabla \vec{R}-(-1)^m \star(\vec{n}\bullet \vec{C}_1 \nabla \log(\rho)) -\nabla f -\nabla^\bot g$ which satisfies 
$$
\left\{\begin{array}{l}
\mbox{div}(\vec{X})=  (-1)^{m} \star( \nabla\vec{n}\bullet \nabla^\bot  \vec{R}) + \nabla^\bot  S\nabla  \star \vec{n} +\mbox{div}( C_0 \nabla \log(r) \star \vec{n}) \\[5mm]
\mbox{curl}(\vec{X})= 0 
\end{array}\right.
$$
Using lemma \ref{aa1}  and lemma \ref{c0est} and lemma \ref{c1e3}, we get that there exist $\vec{\al}$ and $\vec{\beta}$ such that 

\beq
\label{estL21}
\Vert \vec{X}-\vec{\al}\ \nabla \log(\rho) - \vec{\beta}\ \nabla^\bot \log(\rho) \Vert_{L^{2,1}(B(0,1/2)\setminus B(0,2r))} \leq C.\eeq
Proceeding as in (\ref{deltan}),  for all $\rho \in (2r,1/2)$, we get

\beq
\label{fin1}
\left\vert \frac{d\vec{R}_\rho}{d\rho}  -\frac{1}{2\pi} \int_0^{2\pi}  (-1)^m \star(\vec{n}\bullet \frac{\vec{C}_1}{\rho}) \right\vert\leq C \left(\rho \int_{\partial B(0,\rho)} \vert \nabla \vec{n}\vert^2 \; d\sigma \right)^\frac{1}{2} \left(\frac{1}{\rho}\int_{\partial B(0,\rho)}  \vert \nabla \vec{R}\vert^2 +\vert \nabla S \vert^2 \; d\sigma \right)^\frac{1}{2}   +\left\vert \frac{C_0}{\rho}\right\vert
\eeq
In particular, using paragraph \ref{L2} and lemma \ref{c0est}, we get that 
\beq
\label{fin11}
\int_{2r}^{1/2} \left\vert \frac{d \vec{R}_\rho}{d\rho}  -\frac{1}{2\pi} \int_0^{2\pi}  (-1)^m \star\left(\vec{n}\bullet \frac{\vec{C}_1}{\rho}\right)  \right\vert\; d\rho \leq  C ,
\eeq
with $C$ independent of $r$. Then, setting $\vec{f}_\rho=\frac{1}{\vert \partial B(0,\rho)\vert} \int_{\partial B(0,\rho)} \vec{f}\; d\sigma$, since 
$$\frac{d \vec{f}_\rho}{d \rho}= \frac{1}{ 2\pi} \int_0^{2\pi} \frac{\partial \vec{f}}{\partial \rho}(r,\theta) \, d\theta ,$$
we get
\beq
\label{frho}
 \int_{2r}^{1/2} \left\vert  \frac{\partial \vec{f}}{\partial \rho}\right\vert \, d\rho     \leq C \int_{B(0,1)\setminus B(0,r)} \vert \nabla f \vert\, dz \leq C \Vert \nabla f\Vert_2 \leq C.
 \eeq
 Combining (\ref{fin11}) and (\ref{frho}), we have
 \beq
\label{fin11bis}
\int_{2r}^{1/2} \left\vert \frac{d \vec{R}_\rho}{d\rho}  -\frac{1}{2\pi} \int_0^{2\pi}  (-1)^m \star\left(\vec{n}\bullet \frac{\vec{C}_1}{\rho}\right)  -\frac{d f_\rho}{d\rho}\right\vert\; d\rho \leq  C ,
\eeq
Thanks to  (\ref{estL21})  and the convexity of the norm  we get 
\beq
\label{fin2}
\left\Vert \frac{d \vec{R}_\rho}{d\rho}  -\frac{1}{2\pi} \int_0^{2\pi}  (-1)^m \star\left(\vec{n}\bullet \frac{\vec{C}_1}{\rho}\right)-\frac{d f_\rho}{d\rho}  -\frac{\vec{\alpha}}{\rho} \right\Vert_{L^{2,1}(B(0,1/2)\setminus B(0,2r))} \leq  C 
\eeq
Using (\ref{r1}) and the last estimate, we get that 
\beq
\label{fin111}
\begin{split}
&\int_{2r}^{1/2} \left\vert \frac{d  \vec{R}_\rho}{d\rho} -\frac{1}{2\pi} \int_0^{2\pi}  (-1)^m \star\left(\vec{n}\bullet \frac{\vec{C}_1}{\rho}\right) -\frac{\vec{\alpha}}{\rho} \right\vert\; d\rho \\
&\leq \left\Vert \frac{d \vec{R}_\rho}{d\rho}  -\frac{1}{2\pi} \int_0^{2\pi}  (-1)^m \star\left(\vec{n}\bullet \frac{\vec{C}_1}{\rho}\right) -\frac{\vec{\alpha}}{\rho} \right\Vert_{L^{2,1}(B(0,1/2)\setminus B(0,2r))} \left\Vert \frac{1}{\rho} \right\Vert_{L^{2,\infty}(B(0,1/2)\setminus B(0,2r))}   \leq  C,
\end{split}
\eeq
with $C$ independent of $r$. Then using (\ref{fin11bis}) and (\ref{fin111}) , we get
$$\Vert\vec{\alpha}\ \nabla \log(\rho) \Vert_{L^{2,1}(B(0,1/2)\setminus B(0,2r))} \leq C . $$
Hence (\ref{estL21}) gives
\beq
\label{fin22}
\Vert \vec{X} - \vec{\beta}\ \nabla^\bot \log(\rho) \Vert_{L^{2,1}(B(0,1/2)\setminus B(0,2r))} \leq C,
\eeq
with $C$ independent of $r$. Then we set $\vec{X}=\vec{X}_1 \frac{\partial}{\partial \rho} + \vec{X}_2 \frac{\partial}{\rho \partial \theta} $, we get
$$\left\Vert \vec{X}_2 -\frac{\vec{\beta}}{\rho} \right\Vert_{L^{2,1}(B(0,1/2)\setminus B(0,2r))} \leq C, $$ 
with $C$ independent of $r$. Using once more the norm convexity we get from (\ref{estL21}) that  
$$\left\Vert \frac{\vec{\beta}}{\rho} \right\Vert_{L^{2,1}(B(0,1/2)\setminus B(0,2r))} =\left\Vert (\vec{X}_2)_\rho -\frac{\vec{\beta}}{\rho} \right\Vert_{L^{2,1}(B(0,1/2)\setminus B(0,2r))} \leq C,$$
which leads, with (\ref{fin22}), to

$$\Vert \nabla \vec{R}-(-1)^m \star(\vec{n}\bullet \vec{C}_1 \nabla \log(\rho)) -\nabla f -\nabla^\bot g \Vert_{L^{2,1}(B(0,1/2)\setminus B(0,2r))} \leq C ,$$ 
with $C$ independent of $r$. Hence
\beq
\label{ff}
\Vert e^{-\lambda}( \nabla \vec{R}-(-1)^m \star(\vec{n}\bullet \vec{C}_1 \nabla \log(\rho)))\res \nabla ^\bot {\vec{\Phi}} -e^{-\lambda}(\nabla f +\nabla^\bot g)\res \nabla^\bot{\vec{\Phi}}  \Vert_{L^{2,1}(B(0,1/2)\setminus B(0,2r))} \leq C .
\eeq 
To conclude we need to make a last algebraic computation. In order to do so, we  set $$\displaystyle \vec{C}_1 =\sum_{\alpha<\beta} c_{\alpha\beta}\vec{n}_\alpha\wedge \vec{n}_\beta +\sum_\alpha c^1_\alpha \vec{n}_\alpha\wedge  \vec{e}_1 + \sum_\alpha c^2_\alpha \vec{n}_\alpha\wedge \vec{e}_2 +c^3 \vec{e}_1\wedge \vec{e}_2\quad.$$
On the one hand
\beq
\label{ff1}
\langle \vec{C_1}\res \vec{e}_1,\vec{n}_\alpha\rangle=-c^1_\alpha\quad,
\eeq
on the other hand, we get
\be
\begin{split}
\vec{n}\bullet \vec{C}_1 &=\sum_{\alpha<\beta =1}^{m-2} c_{\alpha\beta}((\vec{n}\res \vec{n}_\alpha)\wedge \vec{n}_\beta -(\vec{n}\res \vec{n}_\beta)\wedge \vec{n}_\alpha) +\sum_{\alpha=1}^{m-2} c^1_\alpha (\vec{n}\res \vec{n}_\alpha)\wedge\vec{e}_1 + \sum_{\alpha=1}^{m-2} c^2_\alpha (\vec{n}\res \vec{n}_\alpha)\wedge \vec{e}_2 \\
&=\sum_{\alpha=1}^{m-2} c^1_\alpha \left((-1)^{\alpha-1} \bigwedge_{\beta\not= \alpha} \vec{n}_\beta\right) \wedge\vec{e}_1 + \sum_{\alpha=1}^{m-2}c^2_\alpha \left((-1)^{\alpha-1} \bigwedge_{\beta\not= \alpha} \vec{n}_\beta\right)\wedge \vec{e}_2 \quad.
\end{split}
\ee
Then
$$
\star(\vec{n}\bullet \vec{C}_1) =\sum_{\alpha=1}^{m-2} c^1_\alpha (-1)^{m-1} \vec{e}_2 \wedge \vec{n}_\alpha + \sum_{\alpha=1}^{m-2} c^2_\alpha (-1)^{m} \vec{e}_1 \wedge \vec{n}_\alpha\quad.
$$
Finally we have
\beq
\label{ff2}
\langle \star(\vec{n}\bullet \vec{C}_1) \res \vec{e}_2,\vec{n}_\alpha\rangle =(-1)^{m-1} c^1_\alpha\quad.
\eeq
Hence (\ref{ff1}) and (\ref{ff2}) give
\beq
\label{magic} 
\left\langle e^{-\lambda} (-1)^m \star\left(\vec{n}\bullet \vec{C}_1 \frac{1}{\rho}\right)\res \frac{\partial  {\vec{\Phi}}}{\partial \theta}, \vec{n}^\alpha \right\rangle   = \left \langle e^{-\lambda}  \vec{C}_1 \res  \frac{\partial {\vec{\Phi}} }{\partial \rho},\vec{n}^\alpha\right\rangle  
\eeq
To get our final estimate it suffies to prove that 

$$\Vert e^{\lambda} H^\alpha \Vert_{L^{2}(B(0,1/2)\setminus B(0,2r))} \leq  C\sqrt{ \sup_{\rho\in (r,1)} \left( \int_{B(0,2\rho)\setminus B(0,\frac{\rho}{2})} \vert \nabla \vec{n}\vert^2 \, dz\right)^\frac{1}{2} }$$ 
with $C$ independent of $r$ and $k$. But taking the scalar product of  (\ref{Hfin}) and $\vec{n}^\alpha$, we get 
\beq
\label{ffff}
\left\Vert4 e^{ \lambda} H^\alpha + \left\langle e^{-\lambda} \left(\nabla \vec{R} -  \nabla^\bot \log(\rho) \vec{C}_1\right)\res \nabla^\bot {\vec{\Phi}}, \vec{n}^\alpha \right\rangle \right\Vert_{L^{2,1}(B(0,1/2)\setminus B(0,2r))} \leq C. 
\eeq
Moreover (\ref{ff}) gives
\beq
\left\Vert \left\langle e^{-\lambda}( \nabla \vec{R}-(-1)^m \star(\vec{n}\bullet \vec{C}_1 \nabla \log(\rho))-(\nabla f +\nabla^\bot g))\res \nabla ^\bot {\vec{\Phi}} , \vec{n}^\alpha \right\rangle  \right\Vert_{L^{2,1}(B(0,1/2)\setminus B(0,2r))} \leq C.
\eeq 
Combining (\ref{magic}),the lats inequality gives 
\beq
\label{fff}
\left\Vert \left\langle e^{-\lambda}\left( \nabla \vec{R} + \nabla^\bot \log(\rho) \vec{C}_1)-(\nabla f +\nabla^\bot g)) \right)\res \nabla^\bot {\vec{\Phi}} , \vec{n}^\alpha \right\rangle  \right\Vert_{L^{2,1}(B(0,1/2)\setminus B(0,2r))} \leq C,
\eeq 
with $C$ independent of $r$. Then, combining (\ref{ffff}) and (\ref{fff}) we have that
\beq
\left\Vert4 e^{ \lambda} H^\alpha - e^{-\lambda}\left\langle (2 \vec{C}_1 \nabla^\bot  \log(\rho) -(\nabla f +\nabla^\bot g))\res \nabla^\bot {\vec{\Phi}} , \vec{n}^\alpha \right\rangle  \right\Vert_{L^{2,1}(B(0,1/2)\setminus B(0,2r))} \leq C. 
\eeq
Hence using the duality between  $L^{2,\infty}$ and $L^{2,1}$,  (\ref{fgest}), we see that it suffies to estimate
$$ e^{-\lambda}\left\langle (\vec{C}_1 \nabla^\bot  \log(\rho))\res \nabla^\bot {\vec{\Phi}} , \vec{n}^\alpha \right\rangle. $$
But,
\be
\begin{split}
e^{-\lambda}\left\langle (\vec{C}_1 \nabla^\bot  \log(\rho))\res \nabla^\bot {\vec{\Phi}} , \vec{n}^\alpha \right\rangle &=
\frac{e^{-\lambda}}{\rho} \left \langle \vec{C}_1\res \frac{\partial \vec{\Phi}}{\partial \rho}, \vec{n}^\alpha \right\rangle=
\frac{1}{\rho} \left \langle \vec{C}_1\res \vec{e}_1, \vec{n}^\alpha \right\rangle=\frac{1}{\rho} \left \langle \vec{C}_1, \vec{e}_1 \wedge \vec{n}^\alpha \right\rangle\\
&= \frac{1}{\rho} \left\langle \star\vec{C}_1, \star( \vec{e}_1 \wedge \vec{n}^\alpha)\right\rangle =(-1)^{\alpha+1} \frac{1}{\rho} \left\langle \star\vec{C}_1,  \vec{e}_2 \bigwedge_{\beta\not=\alpha} \vec{n}^\beta\right\rangle \\
&= (-1)^{\alpha+1} \frac{e^{-\lambda}}{\rho} \left\langle \star\vec{C}_1, \frac{\partial \vec{\Phi}}{\rho \partial \theta } \bigwedge_{\beta\not=\alpha} \vec{n}^\beta \right\rangle \\
&= (-1)^{\alpha+1} \frac{e^{-\lambda}}{\rho}  \left( \frac{\partial  \left\langle \star \vec{C}_{1}, \vec{\Phi}\bigwedge_{\beta\not=\alpha} \vec{n}^\beta  \right\rangle}{\rho \partial \theta } -  \left\langle  \star\frac{\partial \vec{C}_1}{\rho \partial \theta },  \vec{\Phi} \bigwedge_{\beta\not=\alpha} \vec{n}^\beta \right\rangle\right.\\
&\left. -\left\langle \star \vec{C}_1,\vec{\Phi} \frac{\partial}{\rho\partial \theta}\left(\bigwedge_{\beta\not=\alpha} \vec{n}^\beta\right) \right\rangle  \right).
\end{split}
\ee
But 
\beq
\label{annul}
\left\langle \star (\vec{c}\wedge \vec{\Phi}),\vec{\Phi}\wedge \left(\bigwedge_{\beta\not=\alpha} \vec{n}^\beta\right) \right\rangle=0.
\eeq
Finally, using (\ref{annul}),  we get 
\be
\begin{split}
e^{-\lambda}\left\langle (\vec{C}_1 \nabla^\bot  \log(\rho))\res \nabla^\bot {\vec{\Phi}} , \vec{n}^\alpha \right\rangle &=
 (-1)^{\alpha+1} \frac{e^{-\lambda}}{\rho}  \left( \frac{\partial  \left\langle \star \vec{c}_{1}, \vec{\Phi}\bigwedge_{\beta\not=\alpha} \vec{n}^\beta  \right\rangle}{\rho \partial \theta } \right.\\
&\left. -  \left\langle  \star\frac{\partial \vec{C}_1}{\rho \partial \theta },  \vec{\Phi} \bigwedge_{\beta\not=\alpha} \vec{n}^\beta \right\rangle  -\left\langle \star \vec{C}_1,\vec{\Phi} \frac{\partial}{\rho\partial \theta}\left(\bigwedge_{\beta\not=\alpha} \vec{n}^\beta\right) \right\rangle\right)\\
&=(-1)^{\alpha+1} \frac{e^{-\lambda}}{\rho}    \left\langle \star \vec{c}_{1},  \frac{\partial\vec{\Phi}}{\rho \partial \theta }\bigwedge_{\beta\not=\alpha} \vec{n}^\beta  \right\rangle + (-1)^{\alpha+1} \frac{e^{-\lambda}}{\rho}    \left\langle \star \vec{c}_{1}, \vec{\Phi} \frac{\partial}{\rho \partial \theta }\bigwedge_{\beta\not=\alpha} \vec{n}^\beta  \right\rangle \\
& - (-1)^{\alpha+1} \frac{e^{-\lambda}}{\rho}  \left\langle  \star\frac{\partial \vec{C}_1}{\rho \partial \theta },  \vec{\Phi} \bigwedge_{\beta\not=\alpha} \vec{n}^\beta \right\rangle  -(-1)^{\alpha+1} \frac{e^{-\lambda}}{\rho} \left\langle \star \vec{C}_1,\vec{\Phi} \frac{\partial}{\rho\partial \theta}\left(\bigwedge_{\beta\not=\alpha} \vec{n}^\beta\right) \right\rangle 
\end{split}
\ee
Moreover, up to a translation we can assume that $\int_{\partial B(0,r)} \vec{\Phi}\, d\sigma  =0$, then thanks to lemma \ref{lambde} and remark \ref{remd}, we easily get that 
$$\left\vert \frac{1}{\vert \partial B(0,\rho) \vert }\int_{\partial B(0,\rho)} \vec{\Phi}\, d\sigma \right\vert \leq C\rho e^{\lambda},$$
Hence $\frac{\vec{\Phi} }{\rho} e^{-\lambda}$  is uniformly bounded.\\

Then, the first term of the right hand side have its $L^2$ norm control by $\lim_{k\rightarrow\infty} \frac{\vec{c}_1^{\,k}}{\sqrt{l_k}}$, the $L^2$-norm control of the second and fourth  terms are controlled by $ \vec{c}_{1} $ and $\Vert \vec{C}_1\Vert_\infty$, which goes to zero when the conformal class degenerates, see section \ref{c1estsec},  Finally the $L^2$-norm of the  third term also goes to zero when the conformal class degenerates thanks to lemma \ref{c1e5}. This achieves the proof of the main theorem.\hfill$\blacksquare$

\section{Proof of theorem~\ref{th-main2}.}

We are exactly under the assumption of the result of Bernard and Rivi\`ere, see \cite{BR}, that we remind for the sake of completeness.\\

\begin{thm}
Let $(\Sigma, h_k)$ a sequence of surfaces with fixed genus, constant curvature and normalized volume if needed. We assume that this sequence converges into the moduli space. Then let  ${\vec{\Phi}}_k : (\Sigma, h_k) \rightarrow \R^m$ a sequence of conformal Willmore immersions with bounded energy, i.e.
$$ \limsup_{k\rightarrow +\infty} W({\vec{\Phi}}_k) <+\infty.$$
Then, there exists a branched smooth immersion ${\vec{\Phi}}_\infty : \Sigma \rightarrow \R^m$ and a finite number of possibly branched  immersions $\omega_j : S^2\rightarrow \R^m$ and $\zeta_t : S^2\rightarrow \R^m$ which are all Willmore  away from possibly finitely many points, and such that, up to a subsequence, 
\begin{equation}
\label{identity}
 \lim_{k\rightarrow +\infty} W({\vec{\Phi}}_k)=  W({\vec{\Phi}}_\infty) + \sum_{j=1}^{p} W(\omega_j) + \sum_{t=1}^{q} (W(\zeta_t)-m_t 4\pi)\quad.
 \end{equation}
where $m_t$ is the integer multiplicity of $\zeta_t$ at the origin.
\hfill $\Box$
\end{thm}

We want to prove that there is no bubble. Let pick a point where a possible bubble concentrates and we choose some conformal parametrization around this point. That is to say, we assume that  $\Phi_k:\D \rightarrow \R^3$ is a conformal Willmore immersion  and concentrates only at $0$. Then by some classical scheme of bubble extraction, we can find $z_k\in \D$, $\lambda_k \in \R_+^*$ and $\theta_0 \in\R$ such that 
$$z_k\rightarrow 0$$
$$\lambda_k \rightarrow 0$$
and 
$$\tilde{\Phi}_k =\frac{\Phi_k(z_k +\lambda_k z)-\Phi_{k}(z_k)}{\vert \Phi_k(z_k) -\Phi_k(z_k+\lambda_ke^{i\theta_0})\vert} \hbox{ converges to } \tilde{\Phi}_\infty \hbox{ in } C^2_{loc}(\R^2) $$
where $\tilde{\Phi}_\infty$ is a conformal Willmore  immersion from $\R^2$ to $\R^3$, which is not a plane.\\

{\bf Claim 1 : $\tilde{\Phi}_\infty$ is bounded at infinity.}\\

If $\tilde{\Phi}_\infty$ is unbounded at infinity, hence considering an appropriate inversion, we can close its image. The closed surface is Lipschitz with a possible branched point, see Lemma A.5 of \cite{Riviere11}, but for topological reason, the order of this branched point is odd, see lemma \ref{topo} in the appendix. If the branching order is bigger than $3$, then the energy will be close to  $12\pi$, hence the closing point is not branched. Moreover the closing point is regular since the residue vanishes (every loop is homotopic to a point), see theorem I.2 in \cite{BR2}. So we have a Bryant sphere whose energy is below $12\pi$ that is to say a round sphere, hence it was a plan before inversion, which is a contradiction.\hfill$\square$\\

{\bf Claim 2 : $\tilde{\Phi}_\infty(\R^2)$ is a round sphere.}\\

Then, $\tilde{\Phi}_\infty$ is bounded at infinity, using  Lemma A.5 of \cite{Riviere11}, it is easy to see that it comes from  a Willmore immersion from $S^2$ to $S^3$ with at most one branched point at infinity.  Hence the first non trivial order is 3, see lemma \ref{topo}, which would insures that this surface has at least $12\pi$ of energy, contradicting our hypothesis. Hence, we can assume there is no branched point at infinity, and it is smooth since the first residue must vanish,  so we have a Bryant whose energy is below $12\pi$ that is to say a round sphere.\hfill$\square$\\

Then, we are going to study the neck region joining this sphere to our limiting surface. Let 
$$l_k (r) = \int _0^{2\pi} \left\vert\frac{\partial \Phi_k}{\partial \theta} \right\vert \,d\theta .$$

Let $r_k$ such that $\displaystyle l(r_k) =\inf_{r\in [\lambda_k, 1/2]} l(r)$.\\

{\bf Claim 3: $\lambda_k << r_k <<1$ for $k$ big enough and $l(r_k)\rightarrow 0$ .}\\

Indeed,  the infimum can't be achieved at the scale $\lambda_k$ of the blown sphere, since we can we can always reduce the length choosing a bigger radius, that it is to say shrinking the curve to the pole. And it is neither achieved at the scale of the limit surface, since  this one is closed once the singularity at $0$ is erased. Moreover $l(r_k)\leq l(r)$ for all $r\in (0,1/2)$ which achieved the proof of the claim.\hfill$\square$\\

Let $\theta_1$ such  that no concentration of the energy occurs at $e^{i\theta_1}$ for $\frac{\Phi_k(z_k +r_k z)}{l(r_k)}$. It is always possible to find such an angle, since there is always only a finite number of  points of concentration. Then,
following our classical scheme of blow-up, we have that 
$$\hat{\Phi}_k =\frac{\Phi_k(z_k +r z)-\Phi_{k}(z_k+e^{i\theta_1}r_k )}{l(r_k)} \hbox{ converges to } \hat {\Phi}_\infty \hbox{ in } C^2_{loc}(\R^2 \setminus S) ,$$
where $S$ is a finite set of possible concentration points. A priori, the conformal factor can diverge. But, since it must satisfies some harnack estimate, if it diverges then it diverges uniformly. If it goes to $-\infty$, we pick $r>0$ such that  there is no concentration on $\partial B(0,r)$, hence the $\Phi_k(\partial B(0,r))$  will be very small with respect to $\Phi_k(\partial B(0,1))$, contradicting the definition of $r_k$. It can't neither go to $+\infty$ since $\Phi_k(\partial B(0,1))$ is of length one.\\

{\bf Claim 4 : $\hat{\Phi}_\infty (\R^2 \setminus S)$ has at least two ends.}\\

We are going to prove that there is one end at $0$, we will obtain the second one at infinity considering an inversion of the parametrization. Indeed, if $\hat{\Phi}_\infty$ where bounded around zero, then we will be able to erase the singularity because the branching order is necessary one and the first residue vanishes. Hence for $r>0$ small enough and $k$ big enough the length of $\Phi_k(\partial B(0,r))$ will be smaller than $1/2$, which will be a contradiction with the definition of $r_k$.\hfill$\square$\\

Now considering an appropriate inversion of $\hat{\Phi}_\infty$, we can close this surface and erase all potential singularities as above. Hence the results is a Willmore sphere which is not a round sphere, else the surface would have only one end before inversion. Thanks to the classification of Bryant, it implies that $\hat{\Phi}_\infty$ have at least four ends.\\

The end corresponding to the infinity is joined to the limiting surface and provides at least $\beta_g$ of energy\footnote{$\beta_g=\inf\{ W(\Phi)\, \vert \, \Phi : \Sigma \rightarrow \R^3 \hbox {immersion and } genus (\Sigma)=g\}$}. The three others are given by $0$ and two other points of concentration. Because they need to be "closed", they contribute to at least $4\pi$ each. It can be seen, thanks to the following monotonicity formula, see \cite{NCRiviere},  let $\Sigma$ a smooth surface with boundary and $\Phi :\Sigma \rightarrow \R^m$ be a weak immersion, then
$$4\pi \leq \int_{\Sigma} \vert \vec{H}\vert^2 \, dv_g + 2 \frac{\mathcal{H}^1(\partial \Sigma)}{d(\partial\Sigma,\Sigma)},$$
where $g=\Phi^*(\xi)$,  $\mathcal{H}^1$ the $1$-dimensional Haussdorff measure, and $\displaystyle d(A,B)=\sup_{p\in A}\inf_{q\in B} \vert p-q\vert+ \sup_{q\in B}\inf_{p\in A} \vert p-q\vert$.\\

Hence applying this monotonicity formula around the three points of concentration that provided an end, we get that total energy of the surface is as close as we want to $\beta_g+ 12\pi$. This excludes the existence of a concentration point and  the proof of the theorem is completed.\hfill $\Box$

\appendix
\section{Residues of Willmore Hopf Tori.}
\label{AHT}
In this section, we briefly explain that the residue $\vec{c}_1$, which notably appears in our main theorem, is a non trivial invariant. Here is the main result
\begin{prop}
\label{AHTP}
Let $\vec{\Phi} :[0,L/2]\times S^1 \rightarrow S^3$ a Willmore Hopf tori, then the residue $\vec{c}_1$ compute on any $\{t\} \times S^1$ vanishes if and only if the tori is equivalent to the  Clifford tori.
\end{prop}
\noindent\noindent{\bf Proof proposition ~\ref{AHTP}.} let $\gamma:[0,L/2] \rightarrow S^2$ an elastic curves such that $\vert \dot{\gamma}\vert=2$, i.e a critical point of
$$ \int_0^{L/2} (1+k^2(s))ds,$$
where $k$ is the geodesic curvature of $\gamma$ into $S^2$.
We consider a lift $\tilde{\gamma}:[0,L/2] \rightarrow S^3$ of $\gamma$, via the Hopf fibration $\pi :S^3\rightarrow S^2$, such that $\tilde{\gamma}$ cuts the fiber of $\pi$ othogonally.
Then we set 
$$\vec{\Phi}(s,\theta)=e^{i\theta} \tilde{\gamma}(s),$$
where $S^3$ is seen as the set of unit quaternions, see \cite{Pinkall}. Then $\vec{\Phi}$ is an isometric immersion, and since $\gamma$ is an elastica, the image of $\vec{\Phi}$ is Willmore in $S^3$, i.e. a critical point of 
$$\int (1+H^2) d\sigma,$$
which implies that $\vec{\Phi}: [0,L/2]\times S^1 \rightarrow \R^4$ is Willmore. In this setting, the normal is equal to
$$\vec{N}=\vec{n}\wedge \vec{\Phi},$$
where $\vec{n}$ is the normal of $\vec{\Phi}$ into $S^3$. We easily check that $\vec{n}=i\partial_s \vec{\Phi}$. The mean curvature vector is equal to
$$\vec{H}=k(s)\vec{n}-\vec{\Phi}.$$
Hence
$$\vec{c}_1= \int_{0}^{2\pi} -(\partial_s \vec{H} -3 \pi_{\vec{N}}(\partial_s \vec{H}) -\star( \partial_\theta \vec{N}\wedge \vec{H}))\wedge \vec{\Phi} +2(\star(\vec{N}\res\vec{H})\res\partial_\theta \vec{\Phi})\, d\theta.$$
Let compute each term separately.
\be
\begin{split}
I&=-\int_{0}^{2\pi} -\partial_s \vec{H}\wedge \vec{\Phi}\, d\theta \\
&=-\int_{0}^{2\pi} \partial_s k \vec{N}-k \partial_s \vec{n} \wedge \vec{\Phi} +\partial_s  \vec{\Phi}\wedge \vec{\Phi} \, d\theta
\end{split}
\ee
we check that 
$$\int_{0}^{2\pi} \vec{N} \, d\theta=  \int_{0}^{2\pi} \vec{n}\wedge \vec{\Phi}  \, d\theta = \pi (\vec{n}(s,0) \wedge \vec{\Phi} (s,0)- \partial_s \vec{\Phi}(s,0) \wedge \partial_\theta   \vec{\Phi}(s,0))$$
and
$$\int_{0}^{2\pi} \partial_s \vec{n} \wedge \vec{\Phi} \, d\theta =-2k \int_{0}^{2\pi} \partial_s \vec{\Phi} \wedge \vec{\Phi} \, d\theta -  \int_{0}^{2\pi} \partial_\theta \vec{\Phi} \wedge \vec{\Phi} \, d\theta $$
where we used the fact that $\partial_s \vec{n}= -2k \partial_s \vec{\Phi}-\partial_\theta \vec{\Phi}$. But
$$\int_{0}^{2\pi} \partial_s \vec{\Phi} \wedge \vec{\Phi} \, d\theta =\pi (\partial_s \vec{\Phi}(s,0) \wedge \vec{\Phi}(s,0)+ \vec{n}(s,0) \wedge \partial_\theta \vec{\Phi}(s,0))$$
and
$$\int_{0}^{2\pi} \partial_\theta \vec{\Phi} \wedge \vec{\Phi} \, d\theta = 2\pi(  \partial_\theta \vec{\Phi}(s,0) \wedge \vec{\Phi}(s,0)) .$$
Finally
\be
\begin{split}
I&= -\pi (\partial_s k (\vec{n}(s,0) \wedge \vec{\Phi} (s,0)- \partial_s \vec{\Phi}(s,0) \wedge \partial_\theta   \vec{\Phi}(s,0))-k(-2k(\partial_s \vec{\Phi}(s,0) \wedge \vec{\Phi}(s,0)+ \vec{n}(s,0) \wedge \partial_\theta \vec{\Phi}(s,0))  \\
&-2(  \partial_\theta \vec{\Phi}(s,0) \wedge \vec{\Phi}(s,0)) )+ (\partial_s \vec{\Phi}(s,0) \wedge \vec{\Phi}(s,0)+ \vec{n}(s,0) \wedge \partial_\theta \vec{\Phi}(s,0)).
\end{split}
\ee
Smilingly we have,
\be
\begin{split}
II&=  \int_{0}^{2\pi} 3 \pi_{\vec{N}}(\partial_s \vec{H})\wedge \vec{\Phi} \, d\theta =3\pi \partial_s k (\vec{n}(s,0) \wedge \vec{\Phi} (s,0)- \partial_s \vec{\Phi}(s,0) \wedge \partial_\theta   \vec{\Phi}(s,0)) \\
&+ 3\pi k(-2k(\partial_s \vec{\Phi}(s,0) \wedge \vec{\Phi}(s,0)+ \vec{n}(s,0) \wedge \partial_\theta \vec{\Phi}(s,0)) -2(  \partial_\theta \vec{\Phi}(s,0) \wedge \vec{\Phi}(s,0)) ) 
\end{split}
\ee

\be
\begin{split}
III=&  \int_{0}^{2\pi}  \star( \partial_\theta \vec{N}\wedge \vec{H}) \wedge \vec{\Phi}  \, d\theta  = \int_{0}^{2\pi} k\partial_\theta \vec{\Phi} \wedge \vec{\Phi} - \partial_s \vec{\Phi} \wedge \vec{\Phi} \, d\theta\\
&=2\pi k(  \partial_\theta \vec{\Phi}(s,0) \wedge \vec{\Phi}(s,0)) -\pi (\partial_s \vec{\Phi}(s,0) \wedge \vec{\Phi}(s,0)+ \vec{n}(s,0) \wedge \partial_\theta \vec{\Phi}(s,0)),
\end{split}
\ee
and
\be
\begin{split}  IV &= \int_{0}^{2\pi} 2(\star(\vec{N}\res\vec{H})\res\partial_\theta \vec{\Phi})\, d\theta = \int_{0}^{2\pi} 2k \partial_\theta \vec{\Phi} \wedge \vec{n} -2 \partial_\theta \vec{\Phi} \wedge \vec{\Phi} \,d\theta \\
&=2\pi k(\partial_\theta \vec{\Phi}(s,0)\wedge \vec{n} (s,0) +\vec{\Phi}(s,0)\wedge \partial_s \vec{\Phi}(s,0))-4\pi(  \partial_\theta \vec{\Phi}(s,0) \wedge \vec{\Phi}(s,0)).
\end{split}
\ee
All the right hand side are orthogonal to $\vec{n}(s,0) \wedge \vec{\Phi} (s,0)$, except the first term of $I$, i.e. $-\pi \partial_s k (\vec{n}(s,0) \wedge \vec{\Phi} (s,0)$. Hence, if $\vec{c}_1=0$, we necessarily get $\partial_s k =0$, but since the equation of an elastica is given by
$$\partial_{ss}k + \frac{1}{2} (k^3+k)=0.$$
We get $k=0$, that is to say $\gamma$ must be a great circle, which implies that $\vec{\Phi}$ parametrizes a Clifford Tori. Reciprocally, it is clear that if $k=0$, then $\vec{c}_1=0$.

\section{Integrability by compensation results.}
First of all, we give a more precise version of lemma A.1 of \cite{LR}.
\begin{lemma}
\label{l1}
Let $0<\eps<\frac{1}{2}$ and $f:\D\setminus B(0,\eps) \rightarrow \R$ an harmonic function.Then, there exists $d\in \R$ and a  positive constant  $C$ independent of $\eps$ and $f$ such that 
$$ \left\Vert \nabla f -\frac{d}{\rho}\right\Vert_{L^{2,1}(B(0,1/2) \setminus B(0,2\eps))} \leq C \Vert \nabla f\Vert_2 .$$
\hfill$\Box$
\end{lemma}
\noindent\noindent{\bf Proof of lemma~\ref{l1}.} We start by decomposing $f$ as a Fourier series, which gives
$$f(\rho,\theta) = c_0 + d_0\ln(\rho) +\sum_{n\in \Z^{*}} (c_n\rho^n+d_n \rho^{-n}) e^{in\theta}.$$
Then we get 
$$f(\rho,\theta)  -c_0 - d_0\ln(\rho)  = \sum_{n\in \Z^{*}}  (c_n\rho^n+d_n \rho^{-n}) e^{in\theta}.$$
Then we estimate the gradient as follows
\be
\left\vert \nabla f(\rho,\theta) -\frac{d_0}{\rho}\right\vert \leq 2 \sum_{n\in \Z^{*}} (\vert n\, c_n\vert \rho^{n-1} + \vert n\, d_n\vert \rho^{-n-1}).
\ee
Then, we estimate  the $L^{2,1}$-norm of the $f_m(z)= \vert z\vert^m$ on $B(0,1/2) \setminus B(0,2\eps)$, for $m\in \Z\setminus \{-1\}$, which gives
\beq
\begin{split}
&\Vert f_m\Vert_{L^{2,1}(B(0,1/2) \setminus B(0,2\eps)} \leq \sqrt{\pi} \int_{0}^{(2\eps)^m} t^{\frac{1}{m}}\; dt \leq 2\sqrt{\pi} (2\eps)^{m+1} \hbox{ for } m<-1\\
&\hbox{ and } \\
&\Vert f_m\Vert_{L^{2,1}(B(0,1/2) \setminus B(0,2\eps)}\leq \frac{\sqrt{\pi}}{2^{m+1}}\hbox{ for } m\geq 0.
\end{split}
\eeq
Hence we get 
\be
\Vert \nabla f \Vert_{L^{2,1}(B(1,2)\setminus B(0,2\eps))} \leq 4\sqrt{\pi}\left( \sum_{n>0} (\vert n\, c_n\vert + \vert n\, d_{-n}\vert )2^{-n}  + \sum_{n<0} ( \vert n\, c_n\vert +\vert n\, d_{-n}\vert)  (2\eps)^{n} \right).
\ee
Hence, thanks to the Cauchy-Scharwz, we get 
\be
\Vert \nabla f \Vert_{L^{2,1}(B(0,1/2) \setminus B(0,2\eps)} \leq 8\sqrt{\pi} \left(\sum_{n\not= 0}  \vert n\vert 2^{-2\vert n\vert} \right)^\frac{1}{2}\left( \sum_{n>0} (\vert n\vert \,\vert  c_n\vert^2  + \vert n\vert \, \vert d_{-n}\vert^2 )  + \sum_{n<0} (\vert n\vert \,\vert  c_n\vert^2  + \vert n\vert \, \vert d_{-n}\vert^2 )   \eps^{2n} \right)^\frac{1}{2} .
\ee
That is to say 
\be
\Vert \nabla f \Vert_{L^{2,1}(B(0,1/2) \setminus B(0,2\eps)} \leq 8\sqrt{\pi} \left(\sum_{n\not= 0}  \vert n\vert 2^{-2\vert n\vert} \right)^\frac{1}{2}\left( \sum_{n>0} (\vert n\vert \,\vert  c_n\vert^2  + \vert n\vert \, \vert d_{n}\vert^2  \eps^{-2n})  + \sum_{n<0} \vert n\vert \,\vert  c_n\vert^2 \eps^{2n}  + \vert n\vert \, \vert d_{n}\vert^2 \right)^\frac{1}{2} .
\ee
Finally we compute the $L^2$-norm of $\nabla f $, 
\be
\begin{split}
\left\Vert \nabla f   \right\Vert_2 &= \vert d_0\vert (\log(1/\eps))^\frac{1}{2} + \left(2\pi \int_{\eps}^{1} \sum_{n\not=0} (\vert n\, c_n\vert^2 \rho^{2n-2} + \vert n\, d_n\vert^2  \rho^{-2n-2} +2 n^2 \Re( c_n \overline{d_n}) \rho^{-2})\,\rho d\rho\right)^\frac{1}{2} \\
&\geq \left(2\pi \sum_{n\not=0}  \frac{\vert n\, c_n\vert^2}{2n} (1-\eps^{2n}) - \frac{\vert n\, d_n\vert^2}{2n}  (1-\eps^{-2n}) +2 n^2 \Re( c_n \overline{d_n}) \log\left(\frac{1}{\eps}\right)\right)^\frac{1}{2} \\
&\geq C \left(\sum_{n<0}  \vert n\vert  \vert c_n\vert^2 \eps^{2n} +  \vert n\vert \vert d_n\vert^2  +\sum_{n>0} \vert n \vert \vert d_n\vert \eps^{-2n} +  \vert n\vert \vert c_n\vert^2 \right)^\frac{1}{2} \\
\end{split}
\ee
which completes the proof of  lemma~\ref{l1} .\hfill$\square$\\

\begin{lemma}
\label{aa1}
Let $0<r<1/2$ and  $a,b,c,d\in W^{1,2}(\D\setminus B(0,r)$ . Let $\vec{X} \in L^1(\D\setminus B(0,r),\R^2)$  satisfying
\be
\left\{\begin{array}{c}\mbox{div}\,(\vec{X}) = a_x\,b_y-a_y \,b_x \\[3mm] \mbox{curl}\,(\vec{X}) = c_x\,d_y-c_y\, d_x\end{array}\right. .
\ee
Then there exists $\vec{\alpha},\vec{\beta} \in \R^2$ and $C>0$ a constant independent of $r$, such that 
$$\Vert \vec{X}-\vec{\alpha} \,\nabla \log(\rho) - \vec{\beta}\, \nabla^\bot \log(\rho) \Vert_{L^{2,1}(B(0,1/2)\setminus B(0,2r))} \leq C(\Vert \nabla a\Vert_2 \Vert \nabla b\Vert_2+\Vert \nabla c\Vert_2\,\Vert \nabla d\Vert_2).$$ 
\hfill $\Box$
\end{lemma}
\noindent{\bf Proof of lemma~\ref{aa1}.} Let $f,g \in W^{1,1}(\D\setminus B(0,r))$ and $\vec{\alpha}',\vec{\beta}' \in \R^2$  such that 
$$\vec{X} =\nabla f + \vec{\alpha}' \ \nabla \log(\rho) +\nabla^\bot g+ \vec{\beta}' \ \nabla^\bot \log(\rho).$$
Then $f$ satisfies
$$\Delta f= a_x\,b_y-a_y \,b_x .$$
Hence we can apply lemma \ref{l1}, and we get that  there exists $\vec{\alpha}  \in \R^2$ and $C>0$ a constant independent of $r$, such that 
$$\Vert \nabla f  -\vec{\alpha} \nabla \log(\rho) \Vert_{L^{2,1}(B(0,1/2)\setminus B(0,2r))} \leq C(\Vert \nabla a\Vert_2 \Vert \nabla b\Vert_2).$$ 
Then we proceed exactly the same way for $g$, which achieves the proof.\hfill$\square$

\medskip

An improvement of the classical Wente inequality was obtain by Bethuel \cite{Bethuel92} using a duality argument.

\medskip

\begin{lemma}
Let $a$ and $b$ be two measurable functions such that $\nabla a \in L^{2,\infty}(B(0,1))$ and $\nabla b \in L^{2}(B(0,1))$. Let $\varphi\in W^{1,1}_0(B_1)$ be the solution of
\be
\Delta\varphi = a_x\,b_y-a_y\,b_x \hbox{ on } B_1
\ee
Then there exists a constant $C$ independent of $\varphi$ such that 
\beq
\label{w3}
\Vert \nabla \varphi \Vert_{2} \leq C\ \Vert \nabla a\Vert_{2,\infty} \Vert \nabla b\Vert_2.
\eeq
\hfill $\Box$
\end{lemma}

\begin{lemma}
\label{lm-a3}
 Let $0<\eps<1/2$, $a$ and $b$ be two measurable functions such that $\nabla a \in L^{2,\infty}(B(0,1) \setminus B(0,\eps) )$ and $\nabla b \in L^{2}(B(0,1)\setminus B(0,\eps) )$. Let $\varphi\in W^{1,1}_0(B(0,1)\setminus B(0,\eps))$ be the solution of
\be
\Delta\varphi = a_x\,b_y-a_y\,b_x \hbox{ on } B_1
\ee
Then there exists a constant $C$ independent of $\varphi$ and $\eps$ such that 
\beq
\label{w3}
\Vert \nabla \varphi \Vert_{2} \leq C \Vert \nabla a\Vert_{2,\infty} \Vert \nabla b\Vert_2.
\eeq
\hfill $\Box$
\end{lemma}
\noindent{\bf Proof of lemma~\ref{lm-a3}.} Let $\tilde{a}$ and $\tilde{b}$ be  Whitney extensions of $a$ and $b$ in $B(0,\eps)$ such that 
$$\Vert  \nabla \tilde{a} \Vert_{2,\infty} \leq C\ \Vert \nabla a\Vert_{2,\infty} $$
and
$$\Vert\nabla \tilde{b} \Vert_2 \leq C\ \Vert \nabla b\Vert_2 . $$
Then, if $\tilde{\varphi} \in  W^{1,1}_0(B(0,1))$ satisfies 
\be
\Delta \tilde{\varphi} = \tilde{a}_x\,\tilde{b}_y-\tilde{a}_y\,\tilde{b}_x \hbox{ on } B(0,1)
\ee
we have from the previous lemma that
\beq
\Vert \nabla \tilde{\varphi} \Vert_{2} \leq C\ \Vert \nabla a\Vert_{2,\infty} \Vert \nabla b\Vert_2.
\eeq
Moreover ${\varphi}-\tilde{\varphi}$ is nothing else than the harmonic extension of $-\,\tilde{\varphi}$ in $B(0,1)\setminus B(0,\eps)$, hence the $L^2$ of its gradient does not exceed the one of $\nabla \tilde{\varphi}$, which permits to conclude.\hfill$\square$
\section{Branched points are of odd order}
We give a short proof of the fact in $\R^3$, branched point obtain as limit of immersion are odd. It can also be found in \cite{LN}, see lemma 3.6.
\begin{lemma}
\label{topo} Let $\vec{\Phi}_k:\D\rightarrow \R^3$ a sequence of immersion which converge to $\vec{\Phi}_\infty:\D\rightarrow \R^3$ in $C^{2}_{loc}(\D\setminus\{0\})$. If $\vec{\Phi}_\infty$ have an isolated branched point of finite order at $0$, then its order is odd.
\end{lemma}

\noindent{\bf Proof of lemma~\ref{topo}.} Up to, Rotate and  reparametrize $\vec{\Phi}_\infty$ we can assume that 

$$\vec{\Phi}_\infty(z)=\left(\begin{array}{c}z^m \\0\end{array}\right)+o(\vert z\vert^{m}) .$$
Let $r>0$  $M_k:\partial B(0,r) \rightarrow SO(3)$ define by
$$M_k(x)= (\vec{e}_1^{\,k},\vec{e}_2^{\,k}, \vec{e}_1^{\,k}\wedge \vec{e}_2^{\,k}),$$
 where $\vec{e}_1^{\,k}= \frac{(\vec{\Phi}_k)_x}{\Vert (\vec{\Phi}_k)_x \Vert }$  and $\vec{e}_2^{\,k}= \frac{(\vec{\Phi}_k)_y- \langle(\vec{\Phi}_k)_y, (\vec{\Phi}_k)_x\rangle }{\Vert (\vec{\Phi}_k)_y- \langle(\vec{\Phi}_k)_y, (\vec{\Phi}_k)_x\rangle \Vert }$.\\
 
Since $\vec{\Phi}_k$ is an immersion, $M_k$ must be in the same connected component than the loop
$$ \Gamma_1(\theta)=\left(\begin{array}{ccc}\cos(\theta) & -\sin(\theta) & 0 \\ \sin(\theta)  & \cos(\theta) & 0 \\0 & 0 & 1\end{array}\right)$$  
But, taking $r>0$ small enough and $k$ large enough, $M_k$ is clearly homotopic to

$$ \Gamma_m(\theta)=\left(\begin{array}{ccc}\cos(m\theta) & -\sin(m\theta) & 0 \\ \sin(m\theta)  & \cos(m\theta) & 0 \\0 & 0 & 1\end{array}\right),$$  
which impose that $m\equiv 1[2]$, since the fundamental group of $SO(3)$ posses exactly two elements: the odd loops and even loops.\hfill$\square$

\end{document}